\documentclass[reqno]{amsart}

\usepackage{amsthm}
\usepackage[usenames,dvipsnames]{pstricks}
\usepackage{epsfig}
\usepackage{pst-grad} 
\usepackage{pst-plot} 
\usepackage{graphicx,subfigure}
\usepackage{amsmath,amssymb}

\usepackage{acronym}
\acrodef{BO}{{\sl Benjamin-Ono}}
\acrodef{rBO}{{\sl regularized Benjamin-Ono}}
\acrodef{rILW}{{\sl regularized Intermediate Long Wave}}
\acrodef{DSW}{{\sl Dispersive Shock Wave}}
\acrodef{DSWs}{{\sl Dispersive Shock Waves}}
\acrodef{ILW}{{\sl Intermediate Long Wave}}
\acrodef{CGN}{{\sl Conjugate Gradient-Newton}}
\acrodef{SW/SW}{{\sl Shallow water / Shallow water}}
\acrodef{B/B}{{\sl Boussinesq / Boussinesq}}

\makeatletter
\newcommand{\sech}{\mathop{\operator@font sech}}
\newcommand{\sign}{\mathop{\operator@font sign}}
\makeatother

\newtheorem{lemma}{Lemma}[section]

\numberwithin{equation}{section}

\begin{document}

\title[]{Solitary-wave solutions of the fractional nonlinear Schr\"{o}dinger equation. II.  A numerical study of the dynamics}


\author{Angel Dur\'an}
\address{\textbf{A.~Dur\'an:} Applied Mathematics Department, University of Valladolid, P/ Belen 15, 47011, Valladolid, Spain}
\email{angel@mac.uva.es}

\author{Nuria Reguera}
\address{\textbf{N.~Reguera:} Department of Mathematics and Computation, University of Burgos, 09001 Burgos, Spain}
\email{nreguera@ubu.es}



\subjclass[2010]{76B25,35C07,65H10}



\keywords{Fractional nonlinear Schr\"{o}dinger equations,  dynamics of solitary waves}

\begin{abstract}
The present paper is a numerical study of the dynamics of solitary wave solutions of the fractional nonlinear Schr\"{o}dinger equation, whose existence was analyzed by the authors in the first part of the project. The computational study will be made from the approximation of the periodic initial-value problem with a fully discrete scheme consisting of a Fourier spectral method for the spatial discretization and a fourth-order, Runge-Kutta-Composition method as time integrator. Several issues regarding the stability of the waves, such as the effects of small and large perturbations, interactions of solitary waves and the resolution of initial data into trains of waves are discussed.
\end{abstract}

\maketitle

\section{Introduction}\label{sec1}
\subsection{The fractional nonlinear Schr\"{o}dinger equation. Solitary wave solutions}
In the first part of the project, \cite{DR1}, the authors considered the 1D version of the focusing fractional nonlinear Schr\"{o}dinger (fNLS) equation
\begin{eqnarray}
iu_{t}- (-\partial_{xx})^{s}u+ |u|^{2\sigma}u=0,\quad {x}\in \mathbb{R},\quad t>0.\label{fnls1d}
\end{eqnarray} 
where $\sigma>0, 0<s<1$. The Fourier multiplier operator $ (-\partial_{xx})^{s}$ has the Fourier symbol
\begin{eqnarray*}
\widehat{(-\partial_{xx})^{s}f}(\xi)=|\xi|^{2s}\widehat{f}(\xi),\quad \xi\in\mathbb{R},\label{fnls2}
\end{eqnarray*}
where
$
\widehat{f}(\xi)$
denotes the Fourier transform of $f\in L^{2}(\mathbb{R})$ at $\xi$. Equation (\ref{fnls1d}) can be written as a real system for $v={\rm Re}u, w={\rm Im}u$
\begin{eqnarray}
v_{t}-(-\partial_{xx})^{s}w+(v^{2}+w^{2})^{\sigma}w&=&0,\nonumber\\
-w_{t}-(-\partial_{xx})^{s}v+(v^{2}+w^{2})^{\sigma}v&=&0.\label{fnls1b}
\end{eqnarray}
Equation (\ref{fnls1d}) is introduced by Laskin in quantum physics, \cite{Laskin2000,Laskin2002,Laskin2011}, by generalizing the Feynman path integrals from stochastic processes of L\'evy motion, and developing a new fractional quantum mechanics. Other applications include the mathematical formulation of Bosom-stars, \cite{FrohlichJL2007}, and some models for the propagation of water waves, \cite{IonescuP2014,ObrechtS2015}. The limiting case $s=1$ corresponds to the classical nonlinear Schr\"{o}dinger (NLS) equation, \cite{SulemS1999}.

Some properties of the initial-value problem (ivp) for (\ref{fnls1d}), related to the purpose of the paper, are here emphasized. The first one is concerned with well-posedness. Among the related literature (see e.~g. \cite{KleinSM2014,DR1}, and references therein), low-regularity well-posedness of (\ref{fnls1d}) for the cubic case ($\sigma=1$) and $1/2<s<1$ is studied in \cite{ChoHKL2015}, for both the ivp and the periodic ivp. The corresponding Cauchy problems are shown to be locally well posed in the Sobolev space $H^{r}$ for $r>s_{g}=\frac{1}{2}(1-s)$, while the nonperiodic problem is ill posed in $H^{r}$ with
$$\frac{1-3s}{2s+1}<r<s_{g}.$$ In \cite{HongS2015}, several results of local well-posedness for the ivp are obtained with respect to $s_{g}$ and $s_{c}=1/2-s/\sigma$. Specifically, local well-posedness in $H^{r}$ holds when $r\geq s_{g}$ and $1/2\leq\sigma<2$ and $r>s_{c}$ when $\sigma\geq 2$ (subcritical cases) and for $r=s_{c}$ where $\sigma>2$ (critical case). The ivp es ill posed in $H^{r}$ for $r\in (s_{c},0)$ when $s\in (1/4,1), \sigma<2s$.

On the other hand, blow-up phenomena for a focusing case like (\ref{fnls1d}) is studied in \cite{Boulenger2016}, where a general criterion for blow-up of radial solutions of the fNLS in $\mathbb{R}^{n}, n \geq 2$ is proved, when $\sigma\geq 2s/n$, being $1/2<s<1$ for $s\geq n/2$ and $2s/n\leq \sigma\leq 2s/(n-2s)$ for $s<n/2$. In the 1D case, when $s\in (1/2,1), s\geq s_{c}>0$, a general blow-up result is proved on a bounded, open interval subject to Dirichlet boundary condition, and from initial conditions with negative energy. Singularity formation in the defocusing 1D case as recently studied, by computational means, in \cite{KleinS}.

In addition (cf. \cite{DR1} and references therein), for smooth, localized solutions, (\ref{fnls1b}) admits a Hamiltonian structure
\begin{eqnarray*}
\frac{\partial}{\partial t}\begin{pmatrix}v\\w\end{pmatrix}=\begin{pmatrix}0&1\\-1&0\end{pmatrix}\delta H(v,w),
\end{eqnarray*}
where $\delta H=(\frac{\delta H}{\delta v},\frac{\delta H}{\delta v})^{T}$ denotes the Fr\'echet derivative and $H$ is the energy function
\begin{eqnarray}
H(v,w)=\int_{\mathbb{R}}\left(\frac{1}{2}\left( (|D|^{s}v)^{2}+(|D|^{s}w)^{2}\right)-\frac{1}{2\sigma+2} (v^{2}+w^{2})^{\sigma+1}\right) d{x},\label{fnls3c}
\end{eqnarray}
with $u=v+iw$ and $|D|^{s}$ is the Fourier multiplier operator satisfying
$$\widehat{|D|^{s}f}(\xi)=|\xi|^{s}\, \widehat{f}(\xi),\quad \xi\in\mathbb{R}.$$
The system (\ref{fnls1b}) admits two other conserved quantities
\begin{eqnarray}
I_{1}(v,w)&=&\frac{1}{2}\int_{\mathbb{R}}(v^{2}+w^{2})dx=\frac{1}{2}\int_{\mathbb{R}}|u|^{2}dx, \label{fnls3a}\\
I_{2}(v,w)&=&\frac{1}{2}\int_{\mathbb{R}}(vw_{x}-wv_{x})dx=\frac{1}{2}\int_{\mathbb{R}}{\rm Im}(u\overline{u}_{x})dx, \label{fnls3b}
\end{eqnarray}
called mass and momentum, respectively, and associated, as in the classical NLS, to the symmetry group, \cite{Olver}, of rotations and translations.

The main purpose of the present paper is the study, by computational means, of the stability and general dynamics of solitary wave solutions of (\ref{fnls1d}). These solitary waves were obtained in \cite{DR1} as critical points of the energy at fixed values of the mass and momentum, of the form
\begin{eqnarray}
\psi(x,t,\lambda_{0}^{1},\lambda_{0}^{2},x_{0},\theta_{0})&=&G_{(t\lambda_{0}^{1},t\lambda_{0}^{2})}(\varphi)\nonumber\\
&=&\rho(x-t\lambda_{0}^{2}-x_{0})e^{i(\theta(x-t\lambda_{0}^{2}-x_{0})+\theta_{0}+\lambda_{0}^{1}t)},\label{fnls22_7}
\end{eqnarray}
where $u_{0}=(v_{0},w_{0})=e^{i\theta(x)}\rho(x)$ with real $\rho$ and $\theta$, and 
$$\varphi(x,x_{0},\theta_{0})=G_{(\theta_{0},x_{0})}(u_{0})=\rho(x-x_{0})e^{i\theta(x-x_{0})+i\theta_{0}},$$ is the {\em orbit} through $u_{0}=(v_{0},w_{0})$ by the symmetry group, and satisfying
\begin{eqnarray}
-(-\partial_{xx})^{s}u_{0}+|u_{0}|^{2\sigma}u_{0}-\lambda_{0}^{1}u_{0}-i\lambda_{0}^{2}\partial_{x}u_{0}=0,\label{fnls22_1}
\end{eqnarray}
or, as a real system
\begin{eqnarray}
-(-\partial_{xx})^{s}v_{0}-\lambda_{0}^{1}v_{0}+\lambda_{0}^{2}w_{0}'+(v_{0}^{2}+w_{0}^{2})^{\sigma}v_{0}&=&0,\nonumber\\
-(-\partial_{xx})^{s}w_{0}-\lambda_{0}^{1}w_{0}-\lambda_{0}^{2}v_{0}'+(v_{0}^{2}+w_{0}^{2})^{\sigma}w_{0}&=&0,\label{fnls22_2}
\end{eqnarray}
with constants (Lagrange multipliers) $\lambda_{0}^{j}, j=1,2$. The choice $\theta(x)=A(x-x_{0})$, for some constant $A$, satisfying
\begin{eqnarray}
\lambda_{0}^{2}=2s|A|^{2s-2}A,\label{aux1}
\end{eqnarray}
$\sigma=1, s\in (1/2,1)$, leads to the subfamily of (\ref{fnls22_7}) considered in \cite{HongS2017}. The existence of smooth solutions of (\ref{fnls22_1}) or (\ref{fnls22_2}), in the general case, was proved in  \cite{DR1}, under the conditions
\begin{eqnarray}
s\in (1/2,1),\quad  \lambda_{0}^{1}>0, \quad |\lambda_{0}^{2}|<c(\lambda_{0}^{1})=2s\left(\frac{\lambda_{0}^{1}}{2s-1}\right)^{\frac{2s-1}{2s}}.\label{fnls_236b}
\end{eqnarray}
Explicit formulas for the waves are not known and accurate numerical procedures for the generation of approximate solitary wave profiles ar required, \cite{DR1}.

The study of several aspects of the dynamics of the solitary wave solutions is the goal of the present paper. The first one concerns the stability under small perturbations. Different types of stability in the literature can be considered. The characterization of the solitary waves as critical points of the Hamiltonian under the constraints of fixed values of mass and momentum invariants can be used to study the stability of the waves by the symmetry group or orbital stability, \cite{GrillakisI,GrillakisII,Weinstein1986,Weinstein1987} (see also the concept of set stability or energetic stability, \cite{Buffoni}). According to the theory developed in \cite{GrillakisII}, the stability depends on the spectrum of the Hessian of the functional
\begin{eqnarray*}
G:=H-\lambda_{0}^{1}I_{1}-\lambda_{0}^{2}I_{2},
\end{eqnarray*}
evaluated at the solitary wave profile $u_{0}=(v_{0},w_{0})$. This leads to
\begin{eqnarray*}
S&:=&G''(v_{0},w_{0})=S_{\infty}+A,\\
S_{\infty}&=&\begin{pmatrix}-(-\partial_{xx})^{s}-\lambda_{0}^{1}&\lambda_{0}^{2}\partial_{x}\\-\lambda_{0}^{2}\partial_{x}&-(-\partial_{xx})^{s}-\lambda_{0}^{1}\end{pmatrix},\\
A&=&\begin{pmatrix}(v_{0}^{2}+w_{0}^{2})^{\sigma-1}(w_{0}^{2}+(2\sigma+1)v_{0}^{2})&2\sigma (v_{0}^{2}+w_{0}^{2})^{\sigma-1}v_{0}w_{0}\\2\sigma (v_{0}^{2}+w_{0}^{2})^{\sigma-1}v_{0}w_{0}&(v_{0}^{2}+w_{0}^{2})^{\sigma-1}(v_{0}^{2}+(2\sigma+1)w_{0}^{2})\end{pmatrix},
\end{eqnarray*}
where we now assume $\sigma\geq 1$. Since $v_{0},w_{0}$ decay to zero at infinity, then $A$ is relatively compact. Therefore, the essential spectrum of $S$ is the same as that of $S_{\infty}$. Fourier analysis is used to compute
\begin{eqnarray*}
\sigma_{\rm ess}(S_{\infty})=\{\lambda=\lambda_{\pm}(\xi)=-(|\xi|^{2s}+\lambda_{0}^{1})\pm\xi\lambda_{0}^{2}, \quad \xi\in\mathbb{R}\}.
\end{eqnarray*}
Therefore, \cite{DR1}, under conditions (\ref{fnls_236b}), it holds that $G''$ is indefinite and the variational theory cannot be applied.

Another type of stability is the asymptotic stability, \cite{MillerW,PegoW,Weinstein1996}. This essentially means that from small enough perturbations of the initial solitary wave profile, the corresponding solutions will evolve asymptotically as
\begin{eqnarray}
\underbrace{\begin{pmatrix}v(x,t)\\w(x,t)\end{pmatrix}}_{u(x,t)}=
\underbrace{\begin{pmatrix}v_{\infty}(x,t,\widetilde{\lambda_{0}^{1}},\widetilde{\lambda_{0}^{2}})\\w_{\infty}(x,t,\widetilde{\lambda_{0}^{1}},\widetilde{\lambda_{0}^{2}})\end{pmatrix}}_{u_{\infty}(x,t,\widetilde{\lambda_{0}^{1}},\widetilde{\lambda_{0}^{2}})}+\underbrace{\begin{pmatrix}z_{1}(x,t)\\z_{2}(x,t)\end{pmatrix}}_{z(x,t)}.\label{as}
\end{eqnarray}
In (\ref{as}), $u_{\infty}$ is a solitary wave associated to modified parameters $\widetilde{\lambda_{0}^{j}}$ close to ${\lambda_{0}^{j}}, j=1,2$, and modified group parameters $\widetilde{x_{0}}, \widetilde{\theta_{0}}$. The second component $z(x,t)$ collects small amplitude, dispersive oscillatory tails as well as smaller nonlinear wave structures. The solitary wave is considered asymptotically stable when the term $z$ tends to zero as $t\rightarrow\infty$ in some sense.

Related to the stability of the solitary waves, some other properties are explored numerically in the present paper, such as the effects of larger perturbations of the solitary waves, interactions of solitary waves and the so-called resolution property, or the resolution of initial data into a series of solitary waves plus decaying small amplitude dispersive tails. This property determines the role of the solitary waves in the general dynamics of the problem, and it has been proved in the case of the soliton solutions of the classical NLS with $\sigma=1$, and other integrable equations, by using the inverse scattering theory, \cite{ZS,BJM}. It has also been studied numerically in the case of other nonlinear dispersive wave problems, cf. e.~g. \cite{DDLM} and references therein.

In this paper, a numerical study of these and other stability issues of the solitary wave solutions  (\ref{fnls22_7}), (\ref{fnls22_1}) is carried out. To this end, the initial-value problem for (\ref{fnls1d}) is first approximated by periodic ivp's on intervals $(-L,L)$ for $L$ long enough. The numerical approximation of the periodic ivp is outlined in Appendix \ref{appA}, and consists of a spectral Galerkin approximation in space and a $4$th-order Runge Kutta Composition method for the time stepping discretization. The fully discrete scheme was already used in \cite{DR1} to check the accuracy of the computation of the approximate solitary wave profiles, and its stability and convergence properties have been shown, either theoretically or numerically, in the approximation of other dispersive nonlinear wave systems, cf. e.~g. \cite{DD2021} and references therein. A numerical analysis of convergence of the method when approximating the periodic ivp for (\ref{fnls1d}) is made in \cite{DR3}.

The paper is structured accordingly to the issues discussed in the numerical study. 
\begin{itemize}
\item The dynamics from small perturbations of different type is analyzed in section \ref{sec2}. The numerical experiments suggest some kind of asymptotic stability of the solitary waves, in the sense, previously described, that the initial perturbed wave evolves asymptotically to a modified waveform plus two groups of oscillatory tails traveling to the right and to the left with respect to the speed of the main wave. The dispersive nature of the tails is theoretically justified from the analysis of plane wave solutions of the corresponding linearized equations, while the nature of the modified wave is discussed, with two possible situations: the asymptotic evolution towards a close solitary wave solution or the approximation to a close breather solution. This type of nwaves has been considered, for the cubic case of (\ref{fnls1d}), in e.~g. \cite{CGLH2019} (see e.~g. \cite{HaragusP} for the cubic focusing NLS).
\item In section 3, the effects of larger perturbations are considered. Here several phenomena can be emphasized. The first one is concerned with a sort of resolution property, which is also observed experimentally in the evolution from different types of initial conditions, like, for example, Gaussian pulses. In some cases, the resolution consists of the formation of a main wave, in the form of a ground state or a moving breather, and solitary wave profiles emerging from it to the right and to the left with respect to this main wave. The generation of some breather type waves is also observed in some experiments involving head-on and overtaking collisions between two solitary waves (which seem inelastic) and between a solitary wave and other waveforms. On the other hand, we checked the possibility of blow-up from large perturbations and under the conditions $\sigma\geq 2s$ and initial data with negative energy, \cite{Boulenger2016}, but this was not observed in any of our experiments, and some phenomena involving moving breather formation was generated instead. 
\item The numerical scheme used to perform the experiments is presented in appendix \ref{appA}, while in appendix \ref{appB} we describe the types of perturbations used in the study (inspired in \cite{DDLM}) as well as some details of the implementation, such as the values of the parameters and the generation of approximate solitary wave profiles. The experiments shown in the paper have used the values $\sigma=1, s=0.8$, but others were made to check that the main features of the dynamics presented here do not change significantly if different values of $s$ and $\sigma$ are considered. As mentioned in appendix \ref{appB}, all the experiments are available from the authors upon request.
\end{itemize}
\subsection{The periodic setting}
As mentioned above, the numerical study of the dynamics of solitary wave solutions of (\ref{fnls1d}) in the present paper is based on the approximation of the ivp with localized initial conditions by a periodic ivp on a long enough interval $(-L,L)$, of the form
\begin{eqnarray}
v_{t}-(-\partial_{xx})^{s}w+(v^{2}+w^{2})^{\sigma}w&=&0,\nonumber\\
-w_{t}-(-\partial_{xx})^{s}v+(v^{2}+w^{2})^{\sigma}v&=&0,\label{fnls1c}
\end{eqnarray}
for $x\in (-L,L), t>0$, with
\begin{eqnarray}
v(x,0)=\widetilde{v}_{0}(x), w(x,0)=\widetilde{w}_{0}(x),\quad x\in (-L,L),\label{fnls1cc}
\end{eqnarray}
smooth, $2L$-periodic given functions. For $T>0$, the ivp (\ref{fnls1c}), (\ref{fnls1cc}) is assumed to be well posed, with a unique, sufficiently smooth solution defined for $t\in [0,T]$. A direct computation proves the following result, concerning the preservation of analogous quantites to (\ref{fnls3c})-(\ref{fnls3b}) for the periodic case:
\begin{lemma}
\label{lemmaNR1}
 The following quantities 
\begin{eqnarray}
I_{1}(v,w)&=&\frac{1}{2}\int_{-L}^{L}(v^{2}+w^{2})dx, \label{fnls3aa}\\
I_{2}(v,w)&=&\frac{1}{2}\int_{-L}^{L}(vw_{x}-wv_{x})dx, \label{fnls3bb}\\
H(v,w)&=&\int_{-L}^{L}\left(\frac{1}{2}\left( (|D|^{s}v)^{2}+(|D|^{s}w)^{2}\right)\right.\nonumber\\
&&\left.-\frac{1}{2\sigma+2} (v^{2}+w^{2})^{\sigma+1}\right) d{x},\label{fnls3cc}
\end{eqnarray}
are invariants by the solution of (\ref{fnls1c}), (\ref{fnls1cc}) for $\widetilde{v}_{0}, \widetilde{w}_{0}$ smooth enough and where  $|D|^{s}$ isthe Fourier multiplier operator  defined by
$$\widehat{|D|^{s}f}(k)=|k|^{s}\, \widehat{f}(k),\quad k\in\mathbb{Z},$$ for $f\in L_{p}^{2}([-L,L])$ squared integrable on $[-L,L]$.
\end{lemma}
The approximation of the ivp by (\ref{fnls1c}), (\ref{fnls1cc})  on long enough intervals $(-L,L)$ has been justified in the literature, cf. e.~g. \cite{Pasciak1982,BChen} and references therein. In this sense, a key property for our purposes here is the asymptotic decay to zero of the solitary wave solutions of (\ref{fnls1d}). This is proved in \cite{DR2}, establishing that, under conditions  (\ref{fnls_236b}), the solution $(v_{0},w_{0})$ of (\ref{fnls22_2}) satisfies
\begin{eqnarray*}
\lim_{|x|\rightarrow\infty}|x|^{2s+1}v_{0}(x)=
\lim_{|x|\rightarrow\infty}|x|^{2s+1}w_{0}(x)=K,
\end{eqnarray*}
for some constant $K$. From the computational point of view, this algebraic decay determines the choice of thr intervals $(-L,L)$ in order to get a good accuracy for long time simulations.
\section{An introductory experiment with the classical NLS}
\label{sec2}
\begin{figure}[htbp]
\centering
\subfigure
{\includegraphics[width=4.1cm]{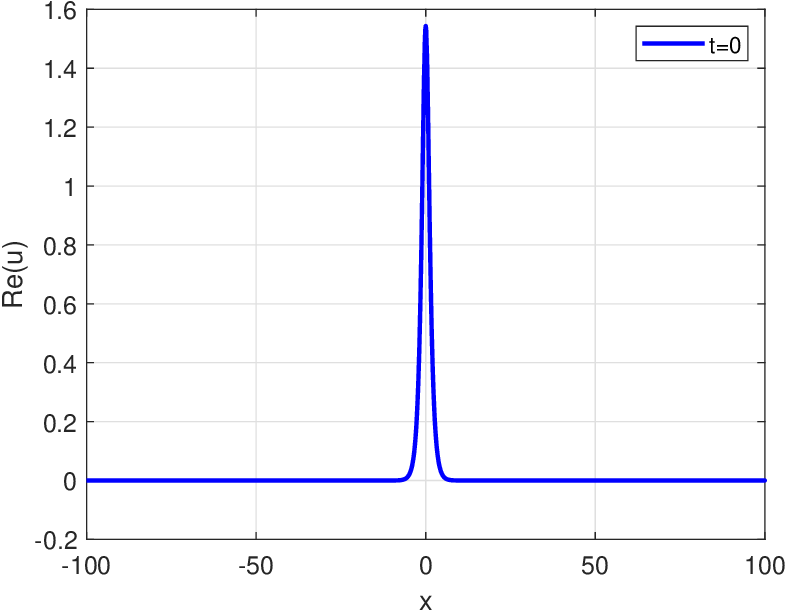}}
\subfigure
{\includegraphics[width=4.1cm]{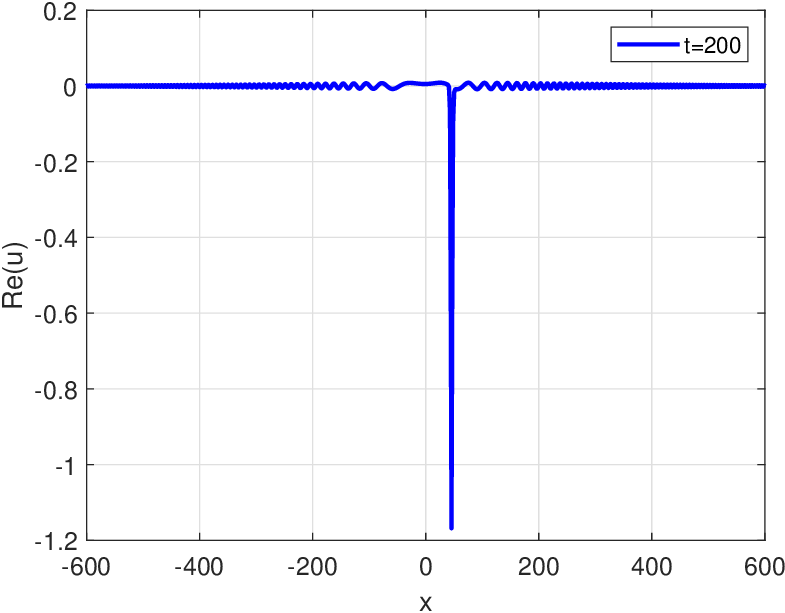}}
\subfigure
{\includegraphics[width=4.1cm]{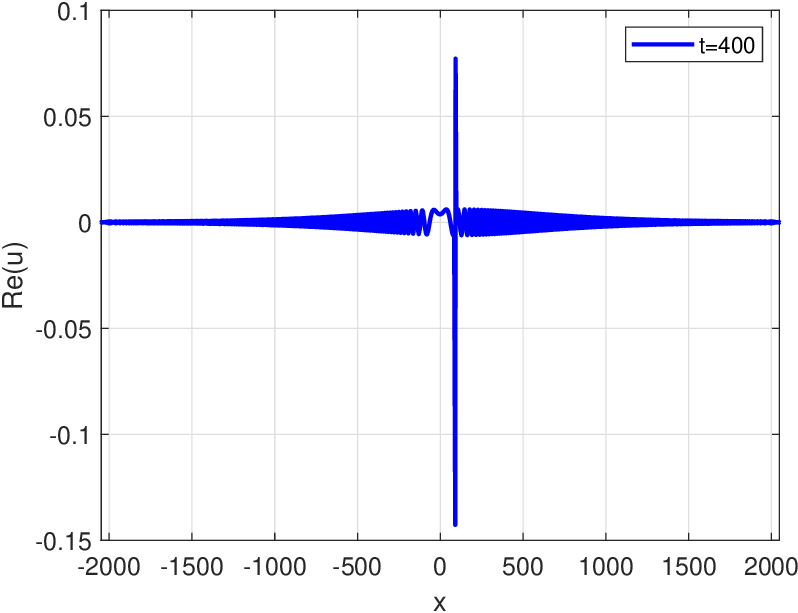}}
\subfigure
{\includegraphics[width=4.1cm]{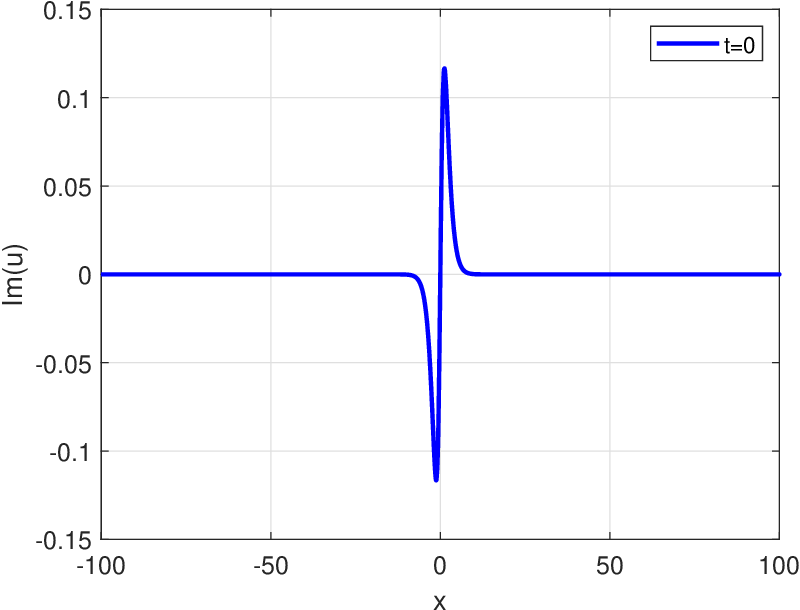}}
\subfigure
{\includegraphics[width=4.1cm]{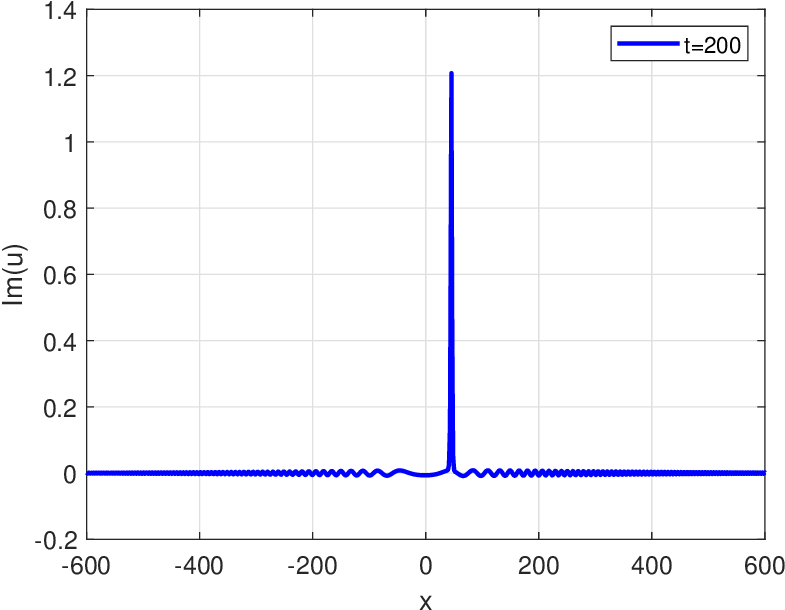}}
\subfigure
{\includegraphics[width=4.1cm]{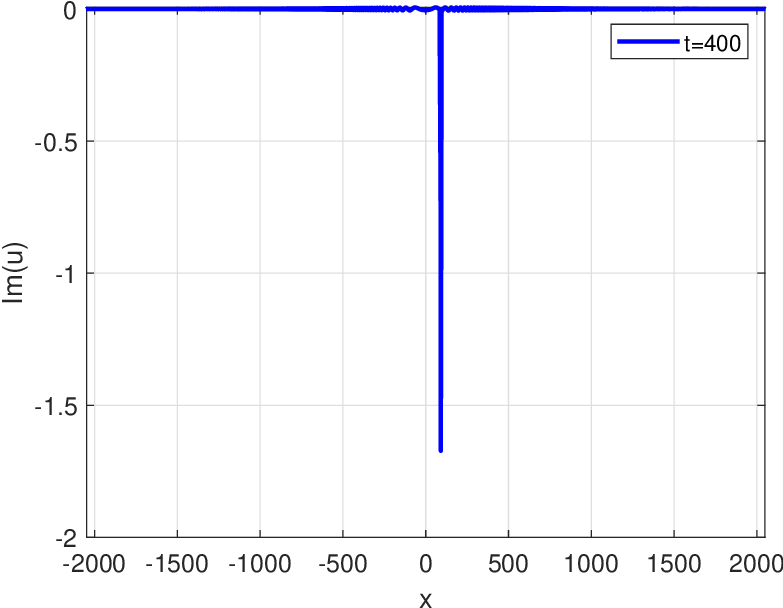}}
\subfigure
{\includegraphics[width=4.1cm]{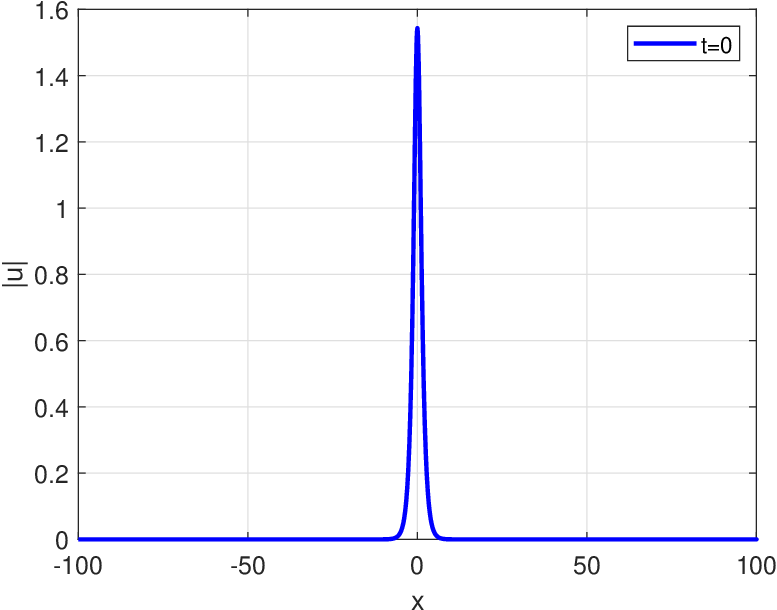}}
\subfigure
{\includegraphics[width=4.1cm]{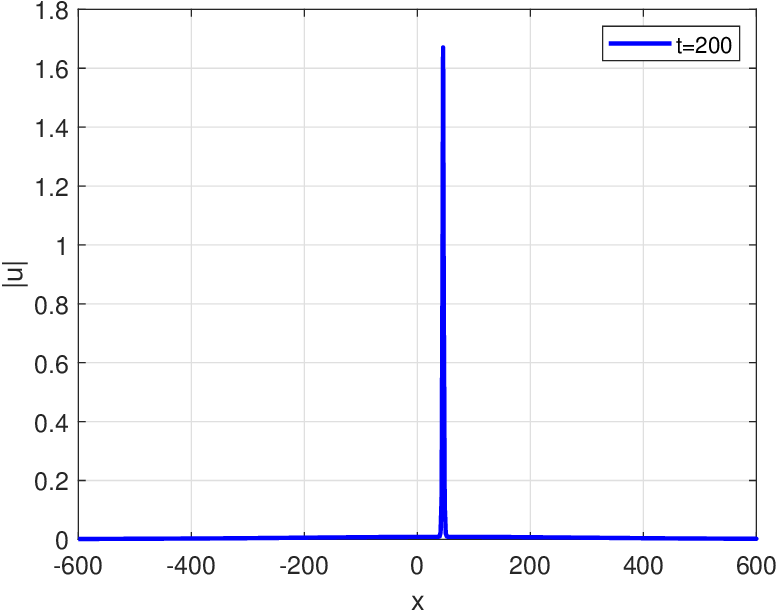}}
\subfigure
{\includegraphics[width=4.1cm]{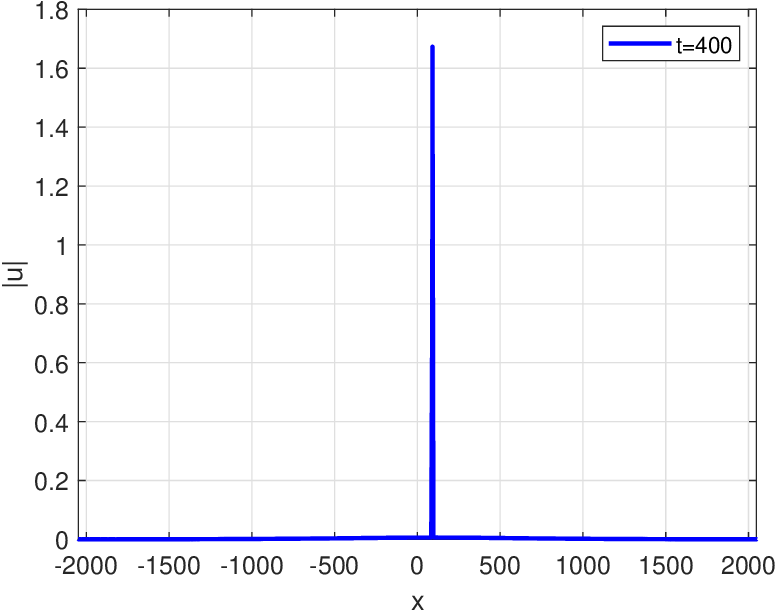}}
\caption{Evolution of the $v, w$, and $\rho$ components of the numerical solution from a slight perturbation (\ref{fnls_b1}) of a solitary wave of the form (\ref{fnls_A5}) in the case $s=\sigma=1$ with $A_{1}=1.1, A_{2}=1$.}
\label{fnls2_FIG0}
\end{figure}

\begin{figure}[htbp]
\centering
\subfigure[]
{\includegraphics[width=6.2cm]{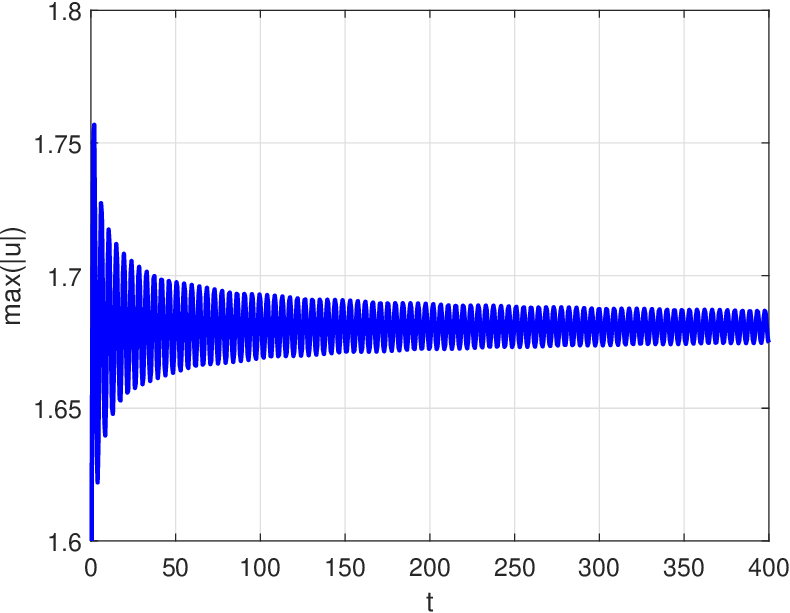}}
\subfigure[]
{\includegraphics[width=6.2cm]{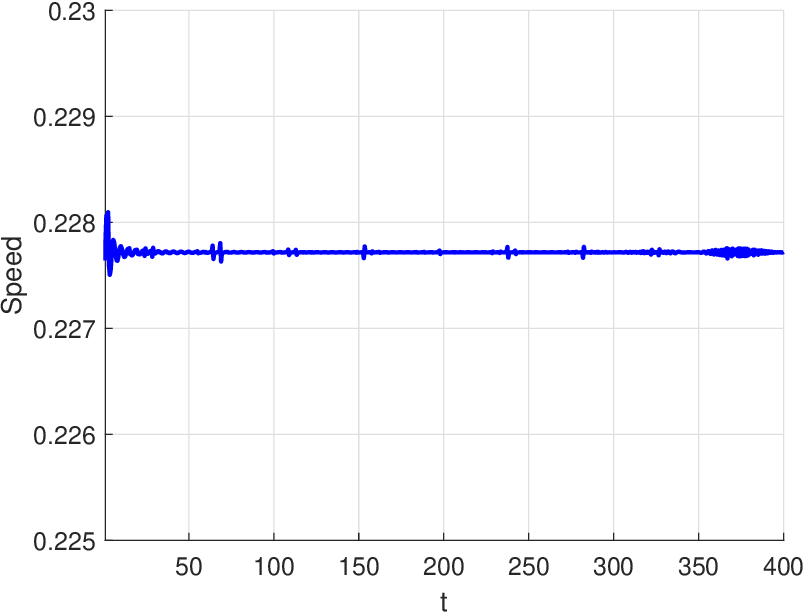}}
\caption{Numerical solution from a slight perturbation (\ref{fnls_b1}) of a solitary wave of the form (\ref{fnls_A5}) in the case $s=\sigma=1$ with $A_{1}=1.1, A_{2}=1$. Emerging wave. Time behaviour of: (a) Amplitude; (b) Speed.}
\label{fnls2_FIG0a}
\end{figure}

By way of comparison, we start the computational study with an experiment concerning the nonfractional, cubic case of (\ref{fnls1d}) (with $s=\sigma=1$) and its soliton solutions. We consider the exact profile (\ref{fnls_A5}) at $t=0$ with $\lambda_{0}^{1}=1, \lambda_{0}^{2}=0.25, x_{0}=\theta_{0}=0$, and perturb it in the form (\ref{fnls_b1}) with $A_{1}=1.1, A_{2}=1$. The perturbed soliton profile is taken as initial condition to run the code described in appendix \ref{appA}. The evolution of the $v, w$ and $\rho$ components of the numerical solution is shown in Figure \ref{fnls2_FIG0}. The figures suggest the formation of a main wave profile traveling to the right, along with dispersive tails (specially observed in the $v$ and $w$ components) leading and trailing the main profile (cf. section \ref{sec32})
The evolution of the amplitude and speed of this wave is displayed in Figure \ref{fnls2_FIG0a}. The behaviour of the amplitude observed in Figure \ref{fnls2_FIG0a}(a) suggests the asymptotic evolution towards a constant value (larger than that of the original profile) while from Figure \ref{fnls2_FIG0a}(b) the wave seems to travel with a slower speed. The soliton solutions (\ref{fnls_A5}) for the cubic case are shown to be asymptotically stable, \cite{CuccagnaP}, in the sense explained in the introduction, and this experiment will serve us to illustrate the asymptotic formation of the modified soliton-type wave and the small-amplitude dispersive tails, and to compare with the subsequent results for the fractional case. We checked that the use of smaller time-step sizes in the numerical integration does not change the  asymptotic oscillatory decrease to some constant amplitude.
\section{Effects of small perturbations}
\label{sec3}
In this section the numerical experiments are focused on the evolution of initial states generated from small perturbations, of different type, of a solitary wave profile. According to the computations, both components evolve to a nearby waveform of different nature, plus two-way propagating dispersive tails. The small size of the perturbation matters, since for larger values other phenomena can be observed, cf. section \ref{sec3}.

\subsection{Evolution of slight perturbations of solitary waves}
\label{sec31}
\begin{figure}[htbp]
\centering
\subfigure
{\includegraphics[width=6.2cm]{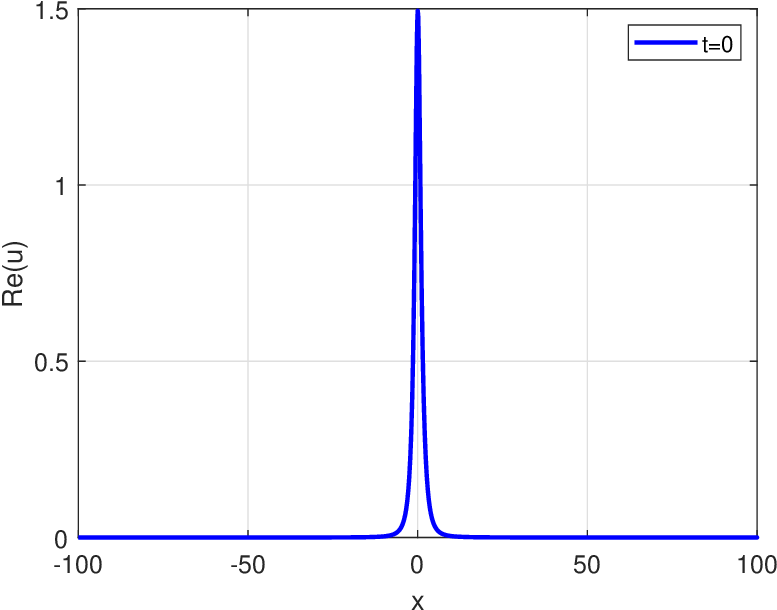}}
\subfigure
{\includegraphics[width=6.2cm]{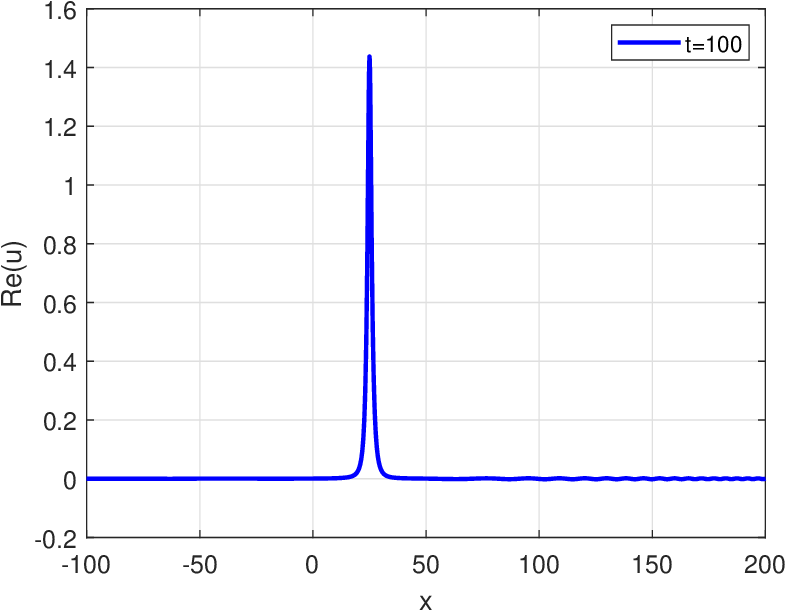}}
\subfigure
{\includegraphics[width=6.2cm]{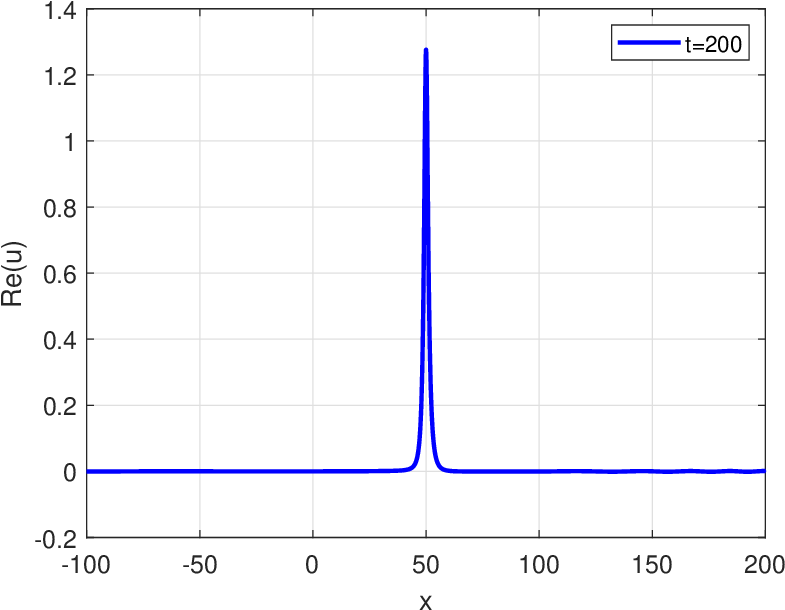}}
\subfigure
{\includegraphics[width=6.2cm]{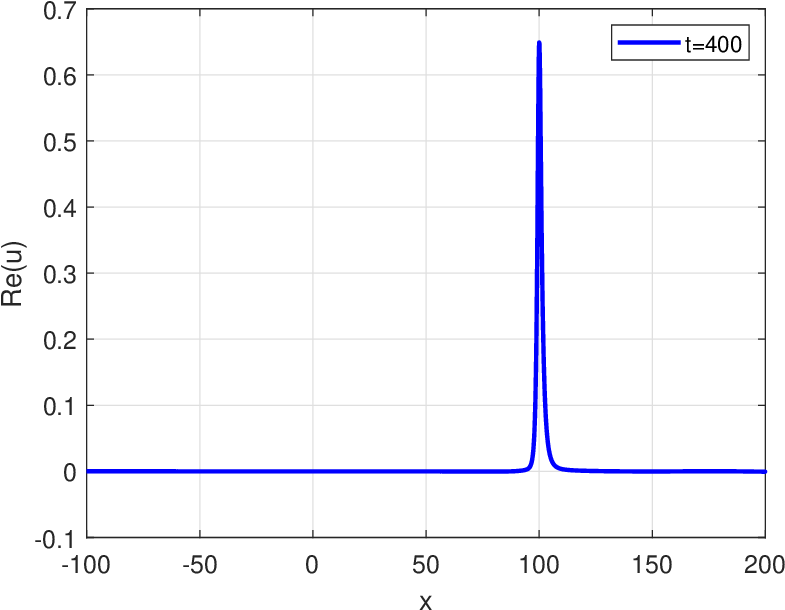}}
\caption{Evolution of the $v$ component of the numerical solution from a slight perturbation of a solitary wave of the form (\ref{fnls_b1}) with $A_{1}=1.1, A_{2}=1$.}
\label{fnls2_FIG1}
\end{figure}

\begin{figure}[htbp]
\centering
\subfigure
{\includegraphics[width=6.2cm]{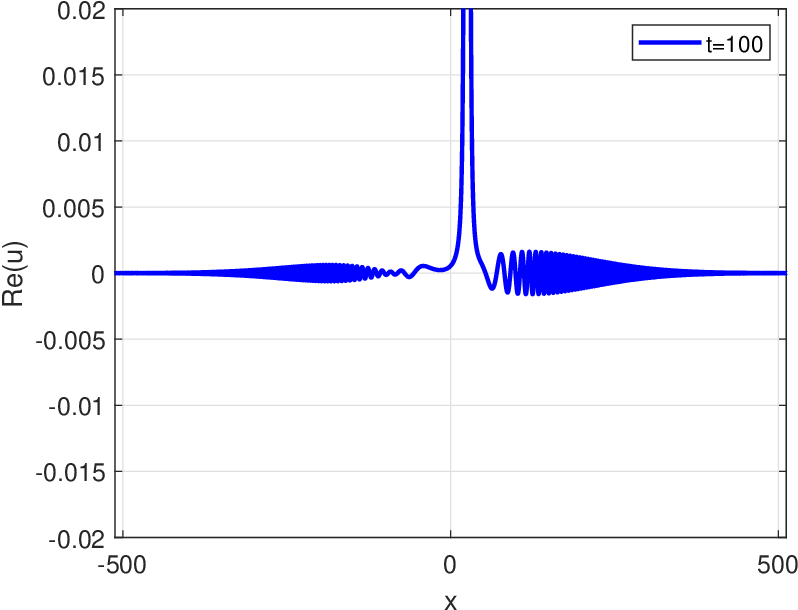}}
\subfigure
{\includegraphics[width=6.2cm]{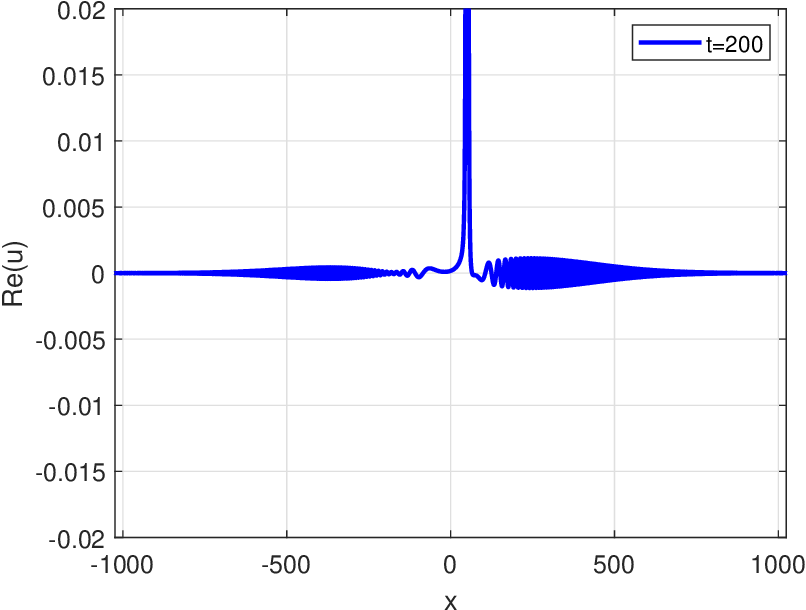}}
\caption{Magnifications of Figure \ref{fnls2_FIG1}.}
\label{fnls2_FIG2}
\end{figure}

\begin{figure}[htbp]
\centering
\subfigure[]
{\includegraphics[width=6.2cm]{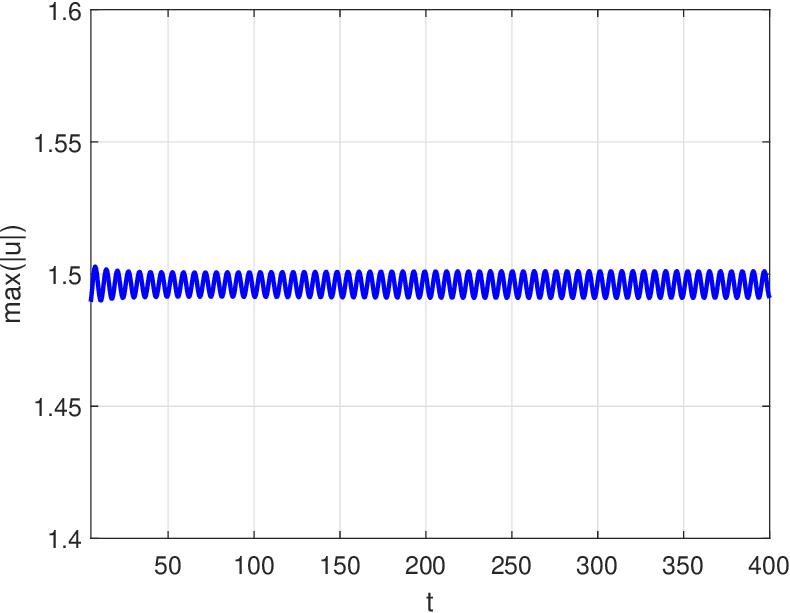}}
\subfigure[]
{\includegraphics[width=6.2cm]{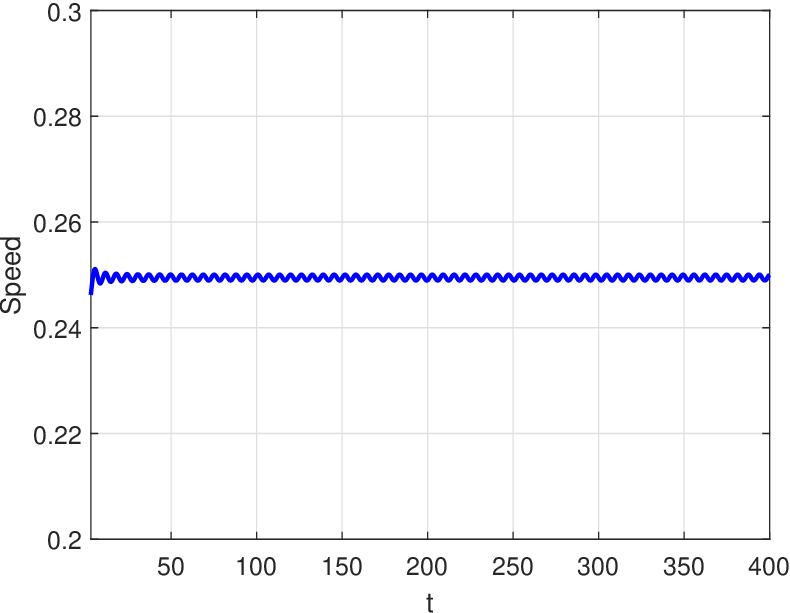}}
\caption{Numerical solution from a slight perturbation of a solitary wave of the form (\ref{fnls_b1}) with $A_{1}=1.1, A_{2}=1$. Emerging wave. Time behaviour of: (a) Amplitude; (b) Speed.}
\label{fnls2_FIG3}
\end{figure}

\begin{figure}[htbp]
\centering
\subfigure
{\includegraphics[width=6.2cm]{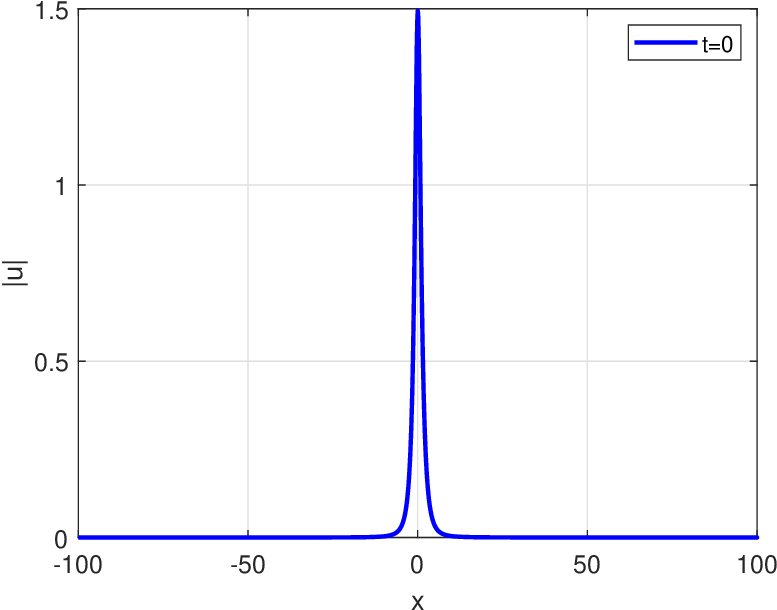}}
\subfigure
{\includegraphics[width=6.2cm]{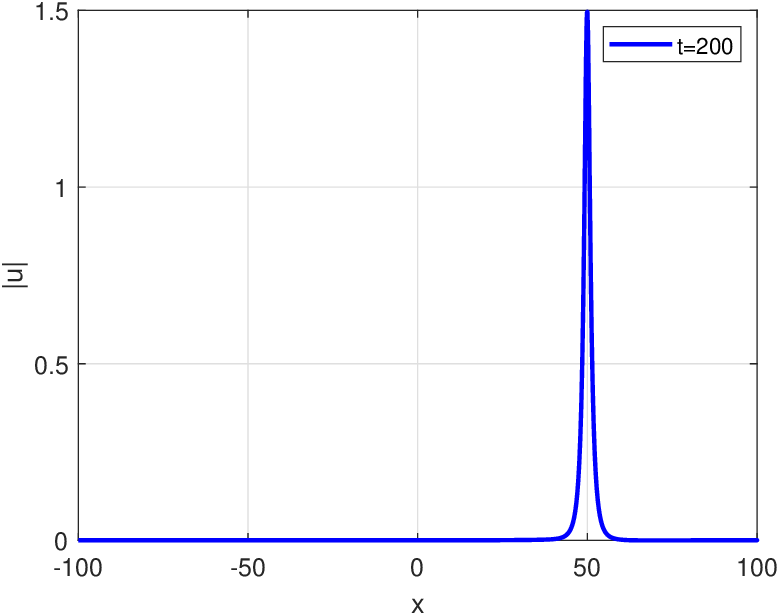}}
\subfigure
{\includegraphics[width=6.2cm]{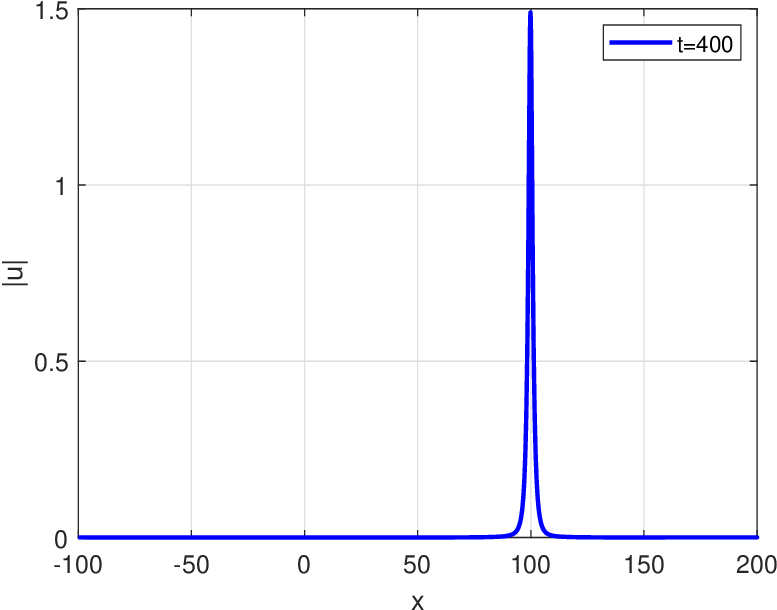}}
\subfigure
{\includegraphics[width=6.2cm]{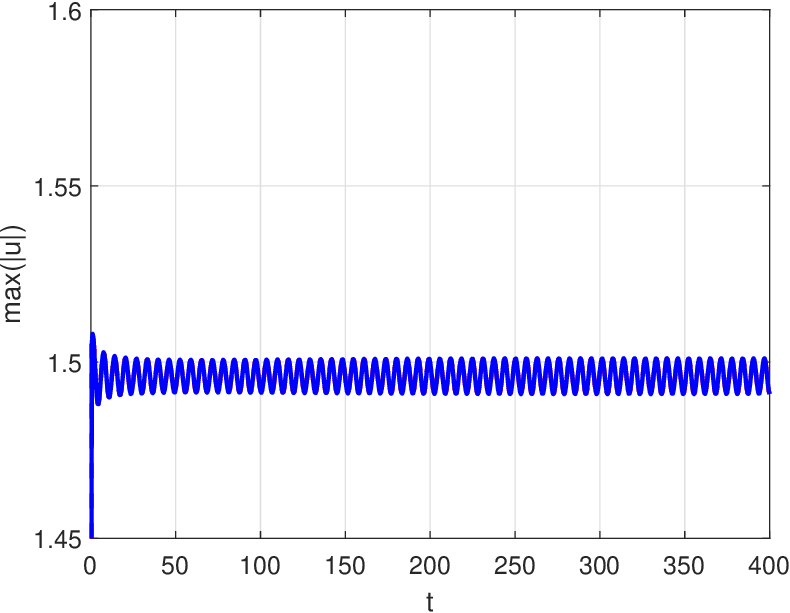}}
\caption{Evolution of the modulus of the numerical solution from a slight perturbation of a solitary wave of the form (\ref{fnls_b1}) with $A_{1}=1.2, A_{2}=1.2$, and time behaviour of the amplitude of the emerging wave.}
\label{fnls2_FIG4}
\end{figure}

A simple experiment, as in the previous section, may consist of generating numerically an approximate solitary wave profile $u_{0}=(v_{0},w_{0})$ and perturbing slightly some of its more relevant parameters, such as the amplitude or the speed, then monitoring the evolution of the numerical solution with the resulting perturbed wave as initial condition. Two examples concerning the amplitude will be shown here. As initial data for the code we take then profiles of the form (\ref{fnls_b1})
with $(A_{1},A_{2})=(1.1,1), (1.2,1.2)$. In the first case, the evolution of the $v$ component of the numerical solution is illustrated in Figure \ref{fnls2_FIG1}, which shows the computed $v$ profile at several times. The unperturbed approximate solitary wave profile is generated with $\lambda_{0}^{2}=0.25$ (playing the role of speed of the wave), $x_{0}=\theta_{0}=0$ and phase $\theta(x)=(x-x_{0})^{2}$. The amplitude of $\rho_{0}=(v_{0}^{2}+w_{0}^{2})^{1/2}$ is $\rho_{0}^{max}=1.4941$. 
The perturbed initial profiles resolves into a single, main waveform, traveling to the right, followed by small amplitude oscillatory tails propagating to the right and to the left compared to that of the main profile, and observed in Figure \ref{fnls2_FIG2}.
The dispersive nature of the emerging tails will be justified theoretically in the following subsection by analizing the plane wave solutions of the corresponding lineatized equations. As far as the nature of the main waveform is concerned, Figure  \ref{fnls2_FIG3} shows the time behaviour of its amplitude and speed. Figure  \ref{fnls2_FIG3}(b) confirms that the wave is moving with a speed close to that of the initial data, while \ref{fnls2_FIG3}(a)  suggests a temporal oscillatory behaviour.
In order to study this behaviour in more detail, Figure  \ref{fnls2_FIG4} illustrates the results obtained with a similar experiment from a perturbed solitary wave profile of the form (\ref{fnls_b1}) with a bit larger perturbation  $(A_{1},A_{2})=(1.2,1.2)$. The evolution of the modulus of the numerical solution confirms the previous dynamics (the oscillatory tails are not shown here) with an amplitude showing some time periodicity, suggesting that the main wave behaves as a moving breather (cf. section \ref{sec2}). The formation and dynamics of breathers from solitons in the cubic case of (\ref{fnls1d}) have been considered in, e.~g. \cite{CGLH2019}.

\begin{figure}[htbp]
\centering
\subfigure[]
{\includegraphics[width=6.2cm]{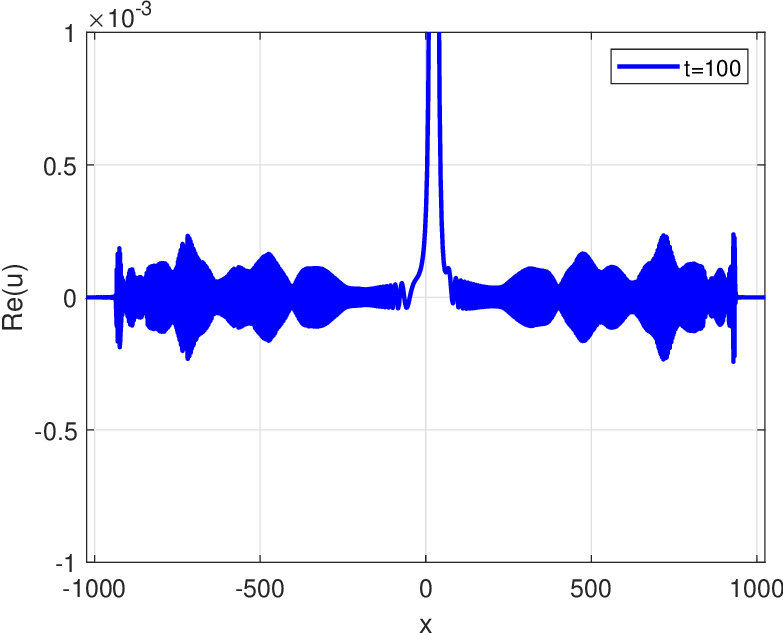}}
\subfigure[]
{\includegraphics[width=6.2cm]{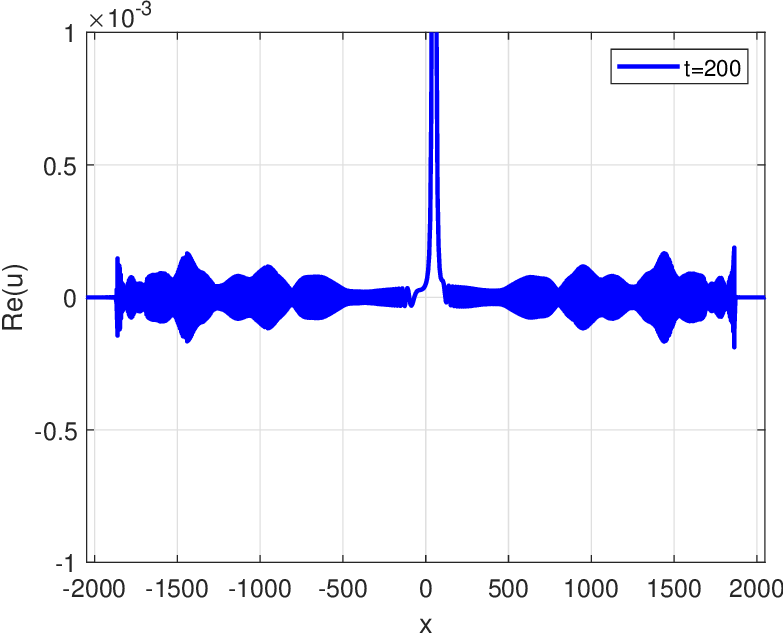}}
\subfigure[]
{\includegraphics[width=6.2cm]{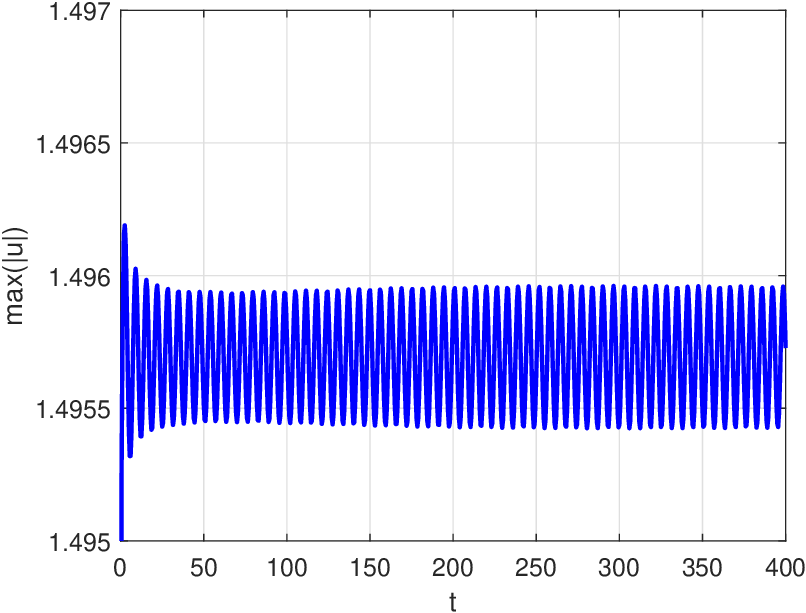}}
\subfigure[]
{\includegraphics[width=6.2cm]{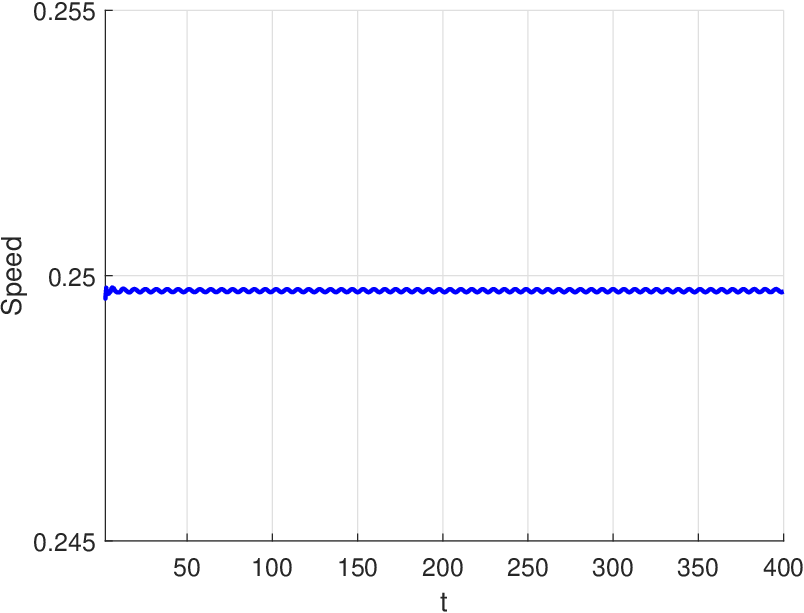}}
\caption{Evolution of the numerical solution from a perturbation of a solitary wave with a numerical noise (\ref{fnls_b3}) where $\beta=10^{5}$. (a), (b) Magnifications of the numerical solution at $t=100,200$; (c), (d) time behaviour of amplitude and speed of the emerging wave.}
\label{fnls2_FIG4a}
\end{figure}

Other numerical experiments with different types of small perturbations of initial solitary wave profiles were performed (cf. Appendix \ref{appB}). By way of illustration, the real part of the computed solitary wave profile has been initially modified with a numerical noise of the form (\ref{fnls_b3}) with $\beta=10^{5}$. Some features of the evolution of the corresponding numerical approximation are shown in Figure \ref{fnls2_FIG4a}. The initial perturbation evolves towards the asymptotic formation of a main wave traveling to the right withb a slower speed (cf. Figure \ref{fnls2_FIG4a}(d)) and time periodic amplitude (cf. Figure \ref{fnls2_FIG4a}(c)). The generation and dynamics of the small-amplitude tails are observed in the magnifications of Figures \ref{fnls2_FIG4a}(a),(b).

The main conclusions concerning the evolution in all cases are similar to those obtained with slight perturbations in amplitude. All the computations suggest a sort of stability, in a sense like that explained in the introduction: the initially perturbed solitary wave evolves to the asymptotic formation of a modified, main waveform (which might be of solitary or moving breather type) along with dispersive, small amplitude tails traveling in both directions (with respect to the main wave).
\subsection{Formation of dispersive tails}
\label{sec32}
Some theoretical arguments to understand the emergence and behaviour of the tails can be given from the idea that small amplitude solutions of the system will approximately evolve according to the linearized equations moving with the solitary wave
\begin{eqnarray*}
(\partial_{t}-c_{s}\partial_{x})v-(-\partial_{yy})^{s}w&=&0,\\
(\partial_{t}-c_{s}\partial_{x})w+(-\partial_{yy})^{s}v&=&0,
\end{eqnarray*}
where $y=x-c_{s}t$, being $c_{s}>0$ the spedd of the wave, which can be simplified to
\begin{eqnarray}
(\partial_{t}-c_{s}\partial_{x})^{2}v+(-\partial_{yy})^{2s}v=0.\label{fnls2_1}
\end{eqnarray}
Plane wave solutions $v(y,t)=e^{i(ky-\omega(k)t)}, k\in\mathbb{R}$, of (\ref{fnls2_1}) will satisfy the dispersion relation
\begin{eqnarray*}
\omega_{\pm}(k)=-c_{s}k\pm\psi(k),\quad \psi(k)=|k|^{2s}, k\in\mathbb{R},
\end{eqnarray*}
for the frequency $\omega(k)$, which can be written as
\begin{eqnarray*}
\omega_{\pm}(k)=-c_{s}k\pm k\varphi(k),\quad \varphi(k)=\left\{\begin{matrix} \frac{|k|^{2s}}{k}&k\neq 0\\0&k=0\end{matrix}.\right.
\end{eqnarray*}
This leads to a local phase speed (relative to $c_{s}$) of the form
\begin{eqnarray*}
V_{\pm}(k)=\frac{\omega_{\pm}(k)}{k}=-c_{s}\pm \varphi(k).
\end{eqnarray*}
Note that $\varphi(k)>0$ (resp. $\varphi(k)<0$) when $k>0$ (resp. $k<0$). In addition, since $2s>1$, then $\varphi(k)$ is increasing with $k>0$ and decreasing for $k<0$. Therefore, for $k>0$
\begin{eqnarray*}
V_{-}(k)<-c_{s}<V_{+}(k),
\end{eqnarray*}
and $V_{+}(k)>0$ for $k>c_{s}^{\frac{1}{2s-1}}$. Thus, part of the components of the plane wave $e^{i(ky-\omega_{+}(k)t)}$ (traveling to the right) leads the main wave and part trails the profile. Similar arguments can be used to show a same property for the waves $e^{i(ky-\omega_{-}(k)t)}$, traveling to the left. This establishes the direction how the plane wave components of the dispersive tails propagate relative to the main wave. Observe that the absolute  phase speed of the wave components traveling to the right $e^{i(ky-\omega_{+}(k)t)}$ ($|V_{+}(k)+c_{s}|$) and to the left $e^{i(ky-\omega_{-}(k)t)}$, ($|V_{-}(k)+c_{s}|$), is not bounded. Components with longer wavelength (smaller $k$) are slower than those of shorter wavelength.

Note also that, since $2s>1$, then $\psi$ is differentiable with
\begin{eqnarray*}
\psi'(k)=\left\{\begin{matrix} 2s k^{2s-1}&k> 0\\-2s k^{2s-1}&k<0\end{matrix}.\right.
\end{eqnarray*}
Then, the associated group velocities are, for $k\neq 0$, given by
\begin{eqnarray*}
\omega_{\pm}'(k)=-c_{s}\pm \psi'(k),
\end{eqnarray*}
with $\psi'(k)\geq 0$ (resp. $\psi'(k)\leq 0$) when $k\geq 0$ (resp. $k\leq 0$), increasing for $k>0$ and decreasing for $k<0$. Thus, for wavenumbers $k>0$, it holds that
\begin{eqnarray*}
\omega_{-}'(k)<c_{s}<\omega_{+}'(k),
\end{eqnarray*}
with $\omega_{+}'(k)>0$ when $k>(c_{s}/2s)^{\frac{1}{2s-1}}$. In the frame of reference $(x,t)$ there are two dispersive groups, one traveling to the left and one to the right (relative to the main wave). In addition, the group velocity is not bounded. All this is illustrated by Figure \ref{fnls2_FIG2}.
\section{Effects of large perturbations}
\label{sec4}
When the size of the perturbations grows, other phenomena can be observed. The first ones are concerned with different types of resolution property, in which solitary waves may emerge, in several ways, during the evolution of large perturbations of solitary waves or other initial data. Furthermore, experiments of interactions between two solitary waves or between a solitary wave and other waveforms may give information on the robustness of these solitary wave structures under more complex perturbations. All these issues will be described in this section.
\subsection{Resolution property}
\label{sec41}
Some computations from different types of initial data suggest a resolution of the initial wave into a sequence of wave sof varied forms plus small-amplitude tails with a main dispersive nature and probably some additional nonlinear structures may emerge from them as times evolves. This resolution property, well understood in the case of integrable equations such as KdV and NLS equations, may be somehow related to the key role of the solitary wave structures in the dynamics of other solutions of the fNLS equation, \cite{Bona}.

\begin{figure}[htbp]
\centering
\subfigure[]
{\includegraphics[width=6.2cm]{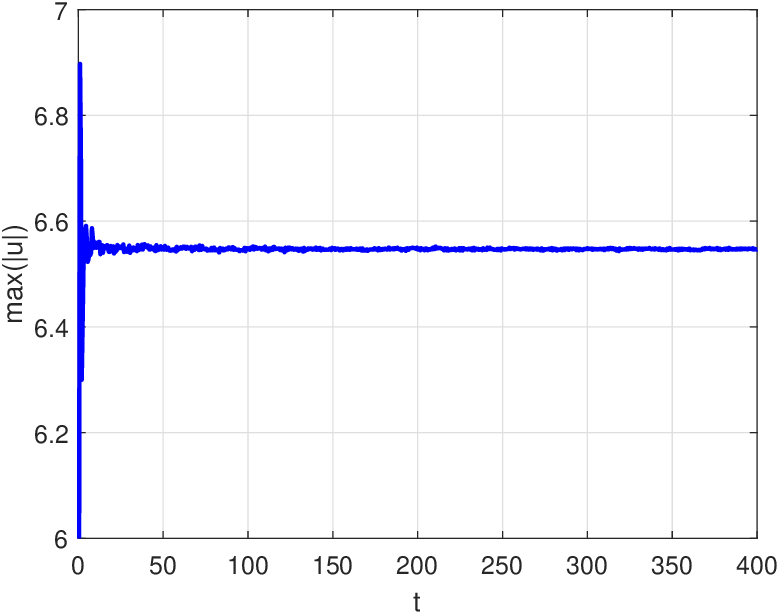}}
\subfigure[]
{\includegraphics[width=6.2cm]{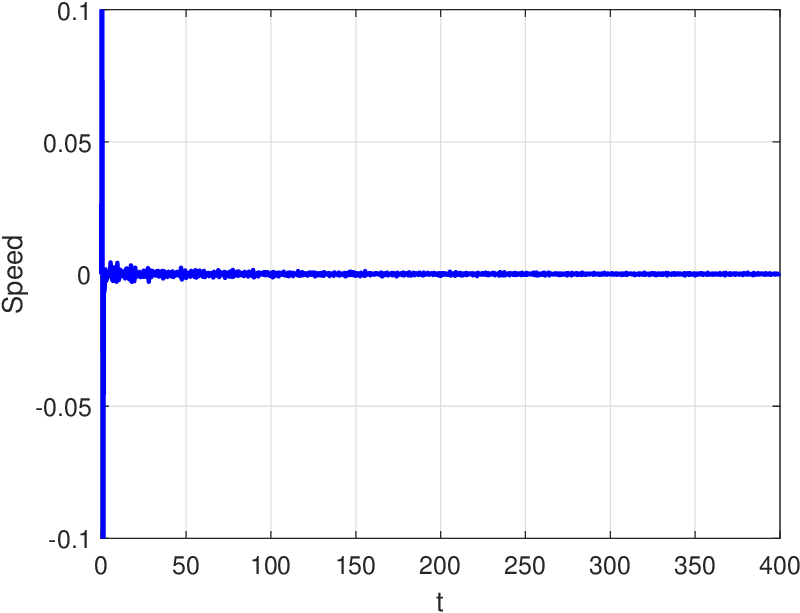}}
\caption{Numerical solution from a slight perturbation of a solitary wave of the form (\ref{fnls_b1}) with $A_{1}=1.8, A_{2}=2$. Tallest emerging wave. Time behaviour of: (a) Amplitude; (b) Speed.}
\label{fnls2_FIG7}
\end{figure}

Two experiments may illustrate the diversity of this phenomenon in this equation. In the first one, we consider a larger perturbation of a solitary wave profile of the form (\ref{fnls_b1}) as initial data, with $A_{1}=1.8, A_{2}=2$ and a linear phase $\theta(x)=A(x-x_{0})$, with $A$ satisfying (\ref{aux1}).
\begin{figure}[htbp]
\centering
\subfigure
{\includegraphics[width=6.2cm]{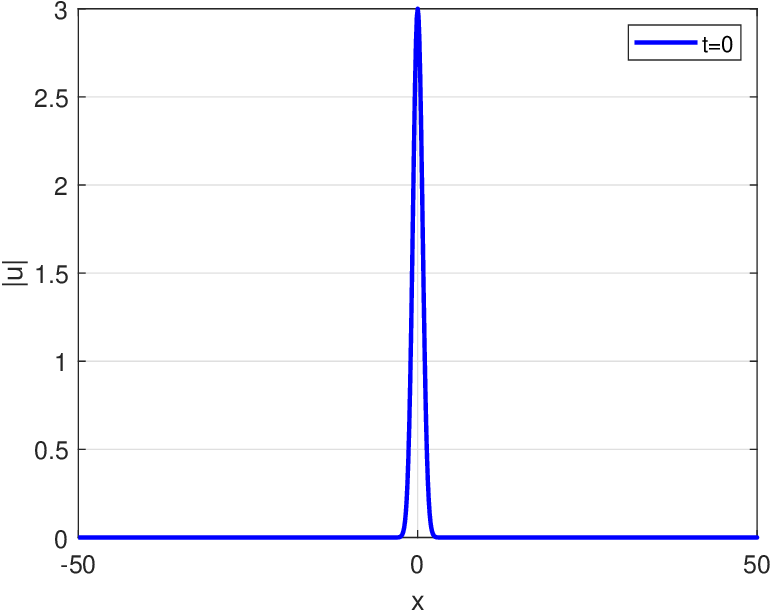}}
\subfigure
{\includegraphics[width=6.2cm]{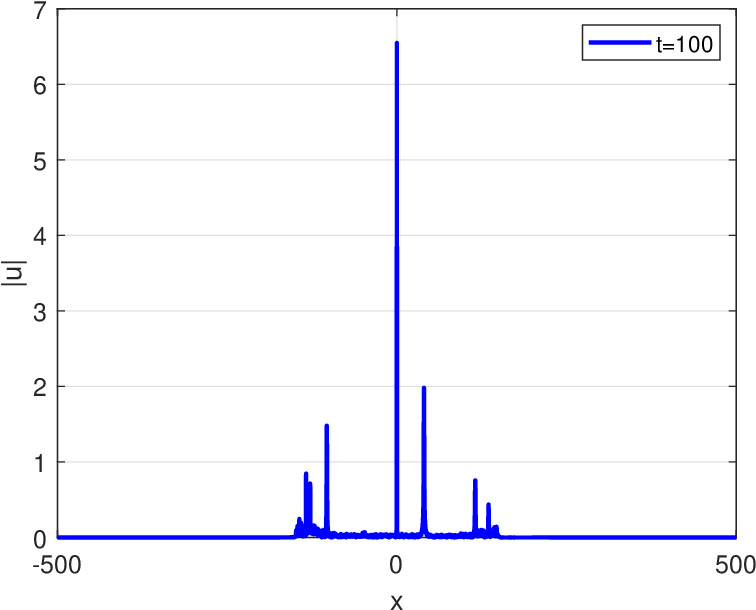}}
\subfigure
{\includegraphics[width=6.2cm]{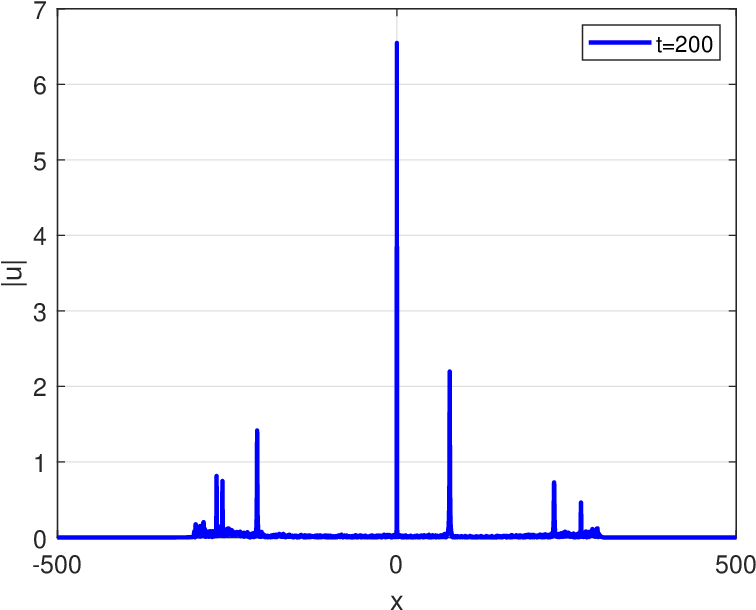}}
\subfigure
{\includegraphics[width=6.2cm]{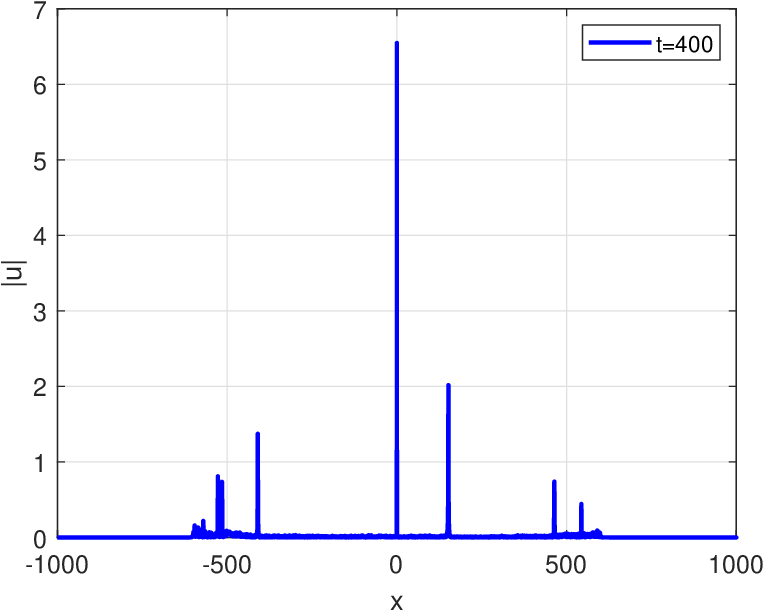}}
\caption{Evolution of the modulus of the numerical solution from a slight perturbation of a solitary wave of the form (\ref{fnls_b1}) with $A_{1}=1.8, A_{2}=2$.}
\label{fnls2_FIG5}
\end{figure}
The evolution of the corresponding approximation is monitored in Figure \ref{fnls2_FIG5}, which displays the modulus of the numerical solution at several times. In this experiment, the resolution consists of the formation of a main wave, with several solitary waves traveling in opposite directions, along with small-amplitude nonlinear waves and dispersive tails, as observed in the magnifications in Figure \ref{fnls2_FIG6}. The behaviour  of the main wave as ground state is suggested in Figure \ref{fnls2_FIG7}, which shows the evolution of its amplitude and speed.

\begin{figure}[htbp]
\centering
\subfigure
{\includegraphics[width=6.2cm]{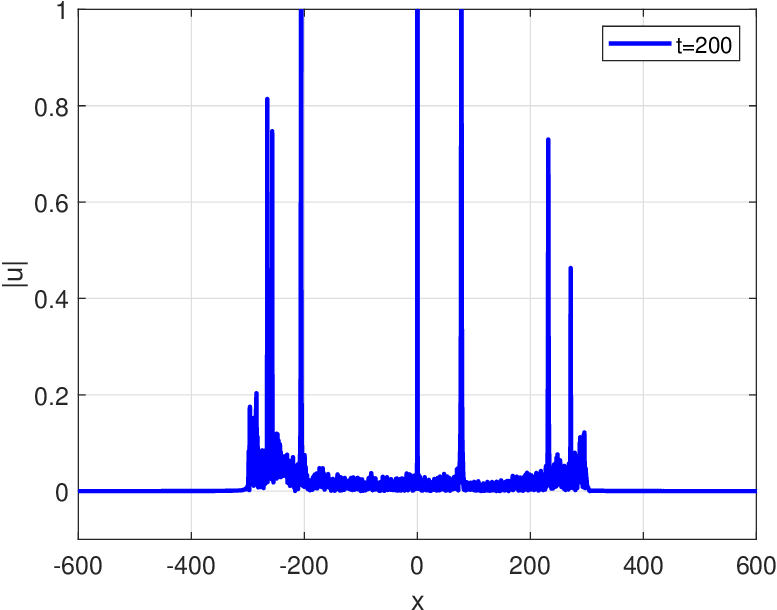}}
\subfigure
{\includegraphics[width=6.2cm]{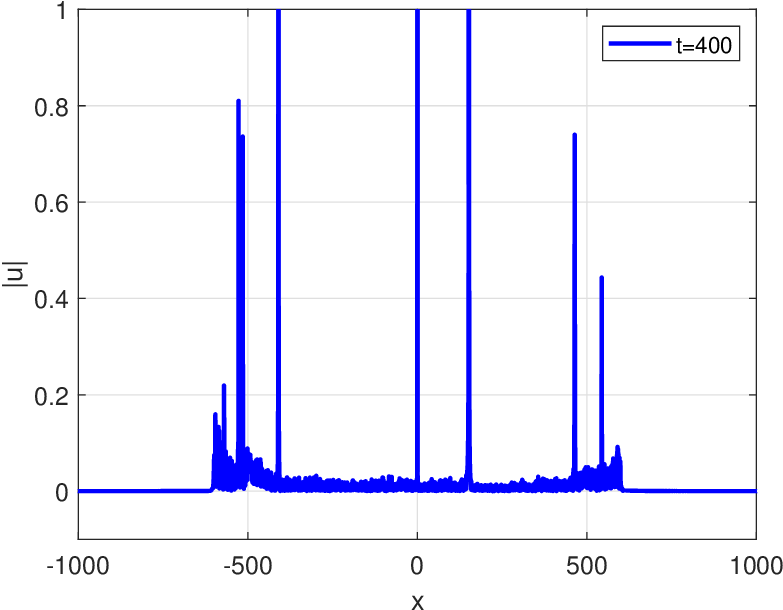}}
\caption{Magnifications of Figure \ref{fnls2_FIG5}.}
\label{fnls2_FIG6}
\end{figure}

\begin{figure}[htbp]
\centering
\subfigure[]
{\includegraphics[width=6.2cm]{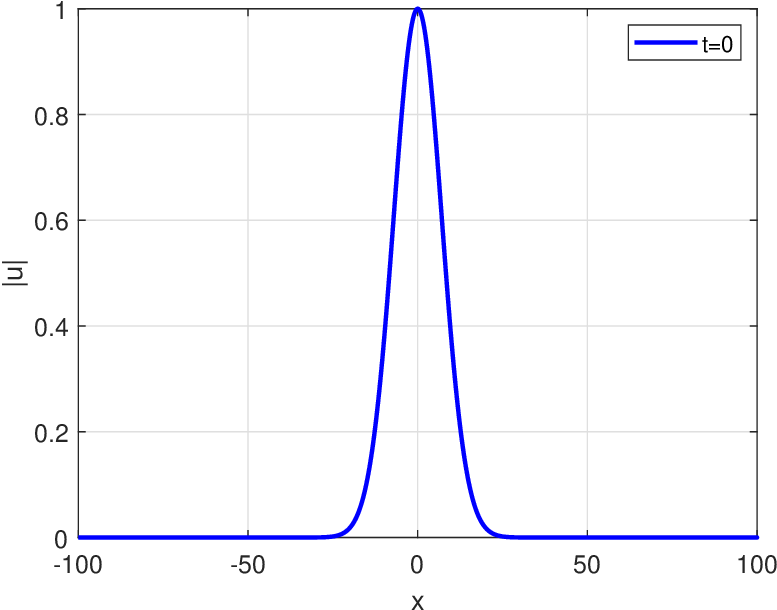}}
\subfigure[]
{\includegraphics[width=6.2cm]{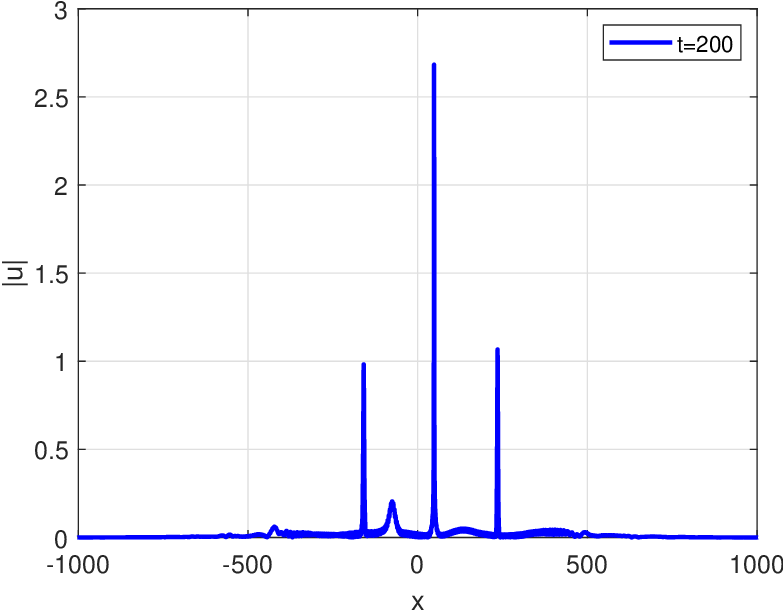}}
\subfigure[]
{\includegraphics[width=6.2cm]{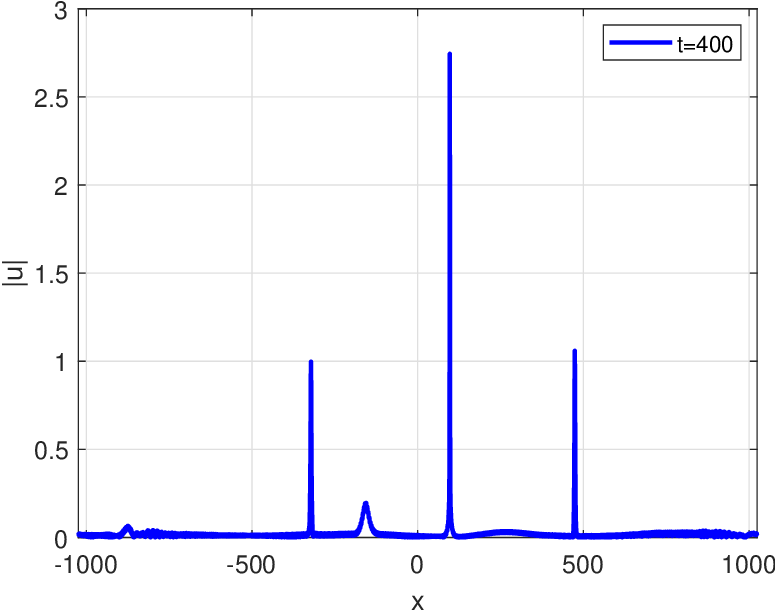}}
\subfigure[]
{\includegraphics[width=6.2cm]{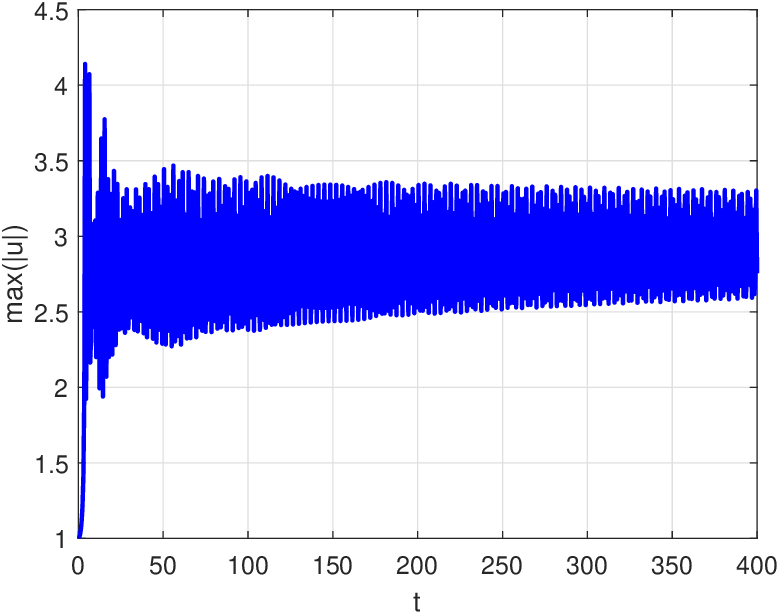}}
\caption{Evolution of the modulus of the numerical solution from an initial condition of the form(\ref{fnls_b5}), with $\theta(x)=A(x-x_{0}), x_{0}=\theta_{0}=0$, and $A_{1}=1, A_{2}=0.01$. (a)-(c) Modulus of the numerical solution at $t=0,200,400$; (d) time behaviour of the amplitude of the tallest emerging wave.}
\label{fnls2_FIG8}
\end{figure}

In the second experiment, we consider an initial condition of the form (\ref{fnls_b5}), with $\theta(x)=A(x-x_{0}), x_{0}=\theta_{0}=0$, and $A_{1}=1, A_{2}=0.01$. The evolution of the modulus of the numercal approximation is given in Figure \ref{fnls2_FIG8}. It is observed that the initial profile breaks into a train of at least two solitary waves traveling to the right and to the left, as well as a main profile moving to the right and with periodic in time amplitude, as suggested in Figure \ref{fnls2_FIG8}(d).
\subsection{Interactions of solitary waves}
\label{sec42}
Another possible manifestation of stability of solitary waves can be suggested from the study of their interactions. It is well known that the nonfractional cubic case ($s=\sigma=1$) is an integrable equation. This implies that solitary wave interactions are elastic, in the sense that after the collision, the waves emerge unchanged in form and speed, compared to the original profiles, with the only modification of a computable phase shift.

\begin{figure}[htbp]
\centering
\subfigure
{\includegraphics[width=6.2cm]{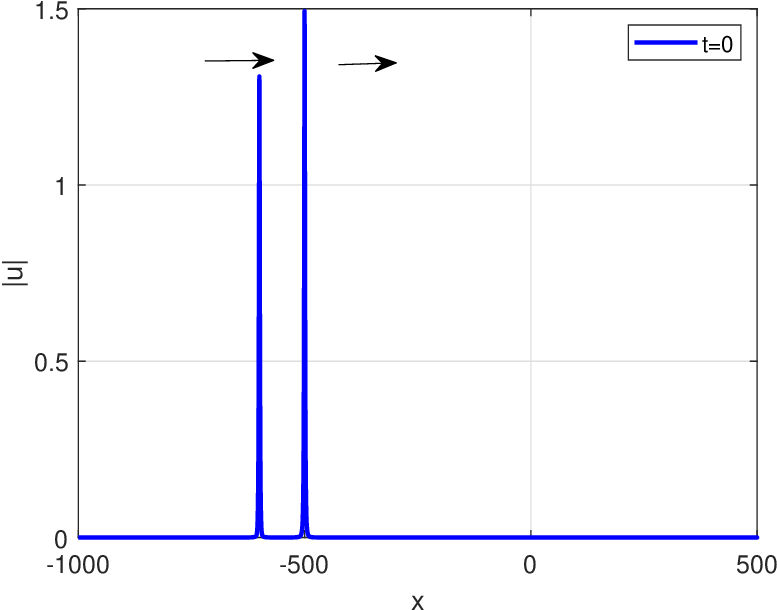}}
\subfigure
{\includegraphics[width=6.2cm]{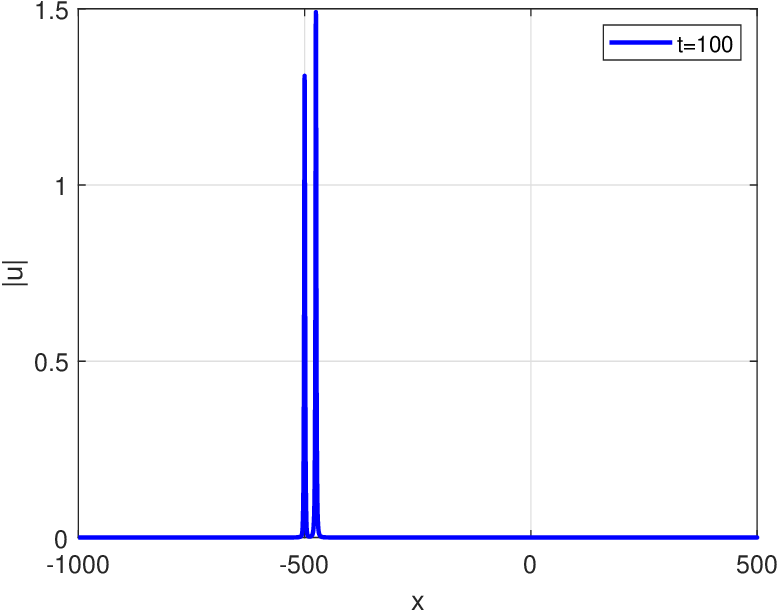}}
\subfigure
{\includegraphics[width=6.2cm]{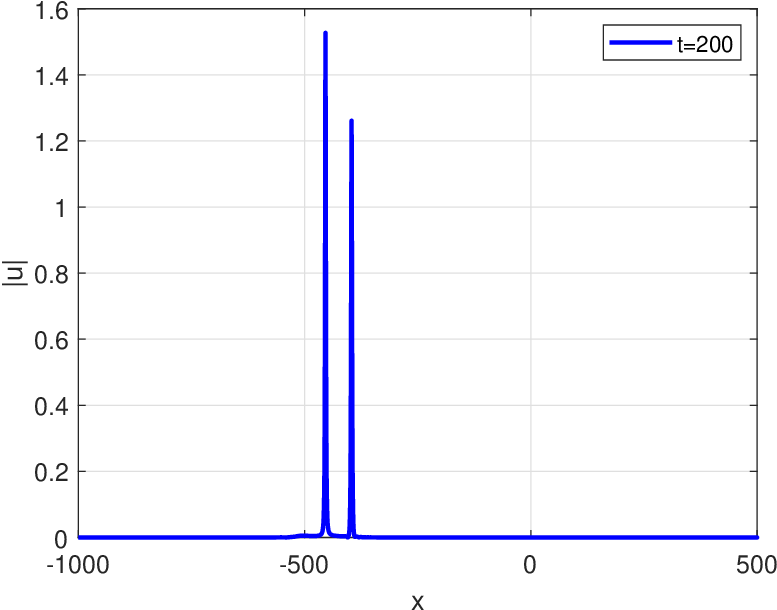}}
\subfigure
{\includegraphics[width=6.2cm]{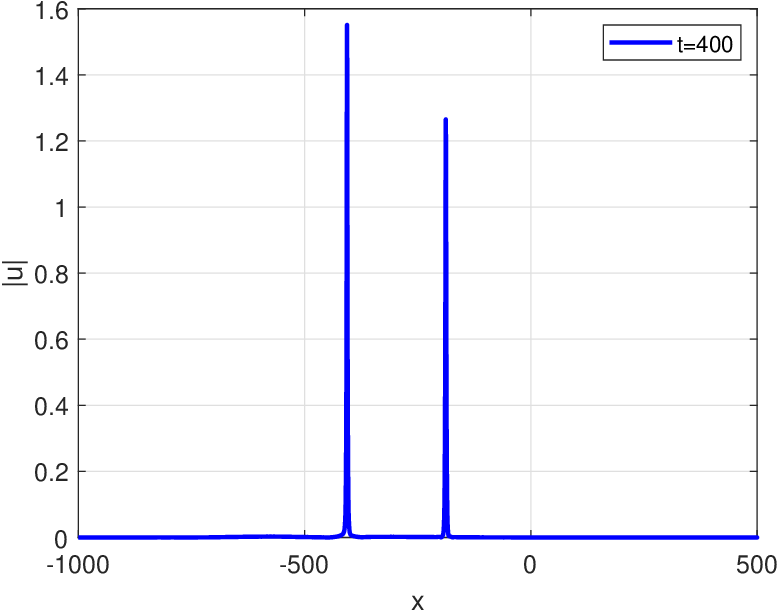}}
\caption{Overtaking collisions. Evolution of the modulus of the numerical solution from a superposition of two computed solitary waves for $\theta(x)=A(x-x_{0})$, $a$ satisfying (\ref{aux1}), with $\lambda_{0}^{1}=1$ and $\lambda_{0}^{2}=1, 0.25$, centered at $x_{0}=-600,-500$ respectively.}
\label{fnls2_FIG9}
\end{figure}
\begin{figure}[htbp]
\centering
\subfigure[]
{\includegraphics[width=6.2cm]{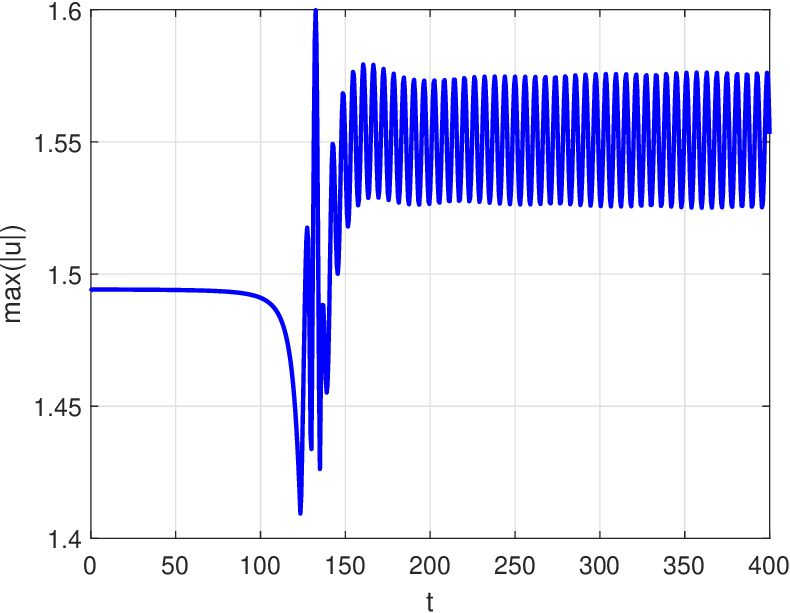}}
\subfigure[]
{\includegraphics[width=6.2cm]{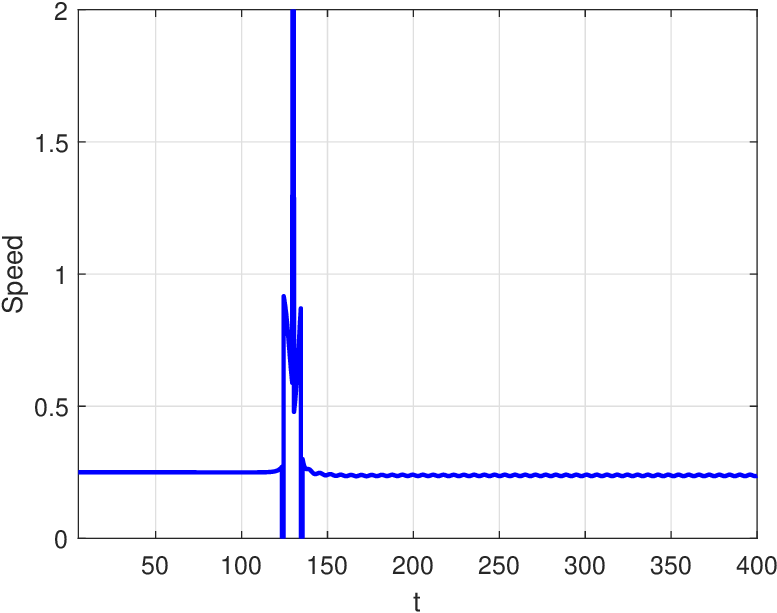}}
\caption{Overtaking collisions. Numerical solution from a superposition of two computed solitary waves for $\theta(x)=A(x-x_{0})$, $a$ satisfying (\ref{aux1}), with $\lambda_{0}^{1}=1$ and $\lambda_{0}^{2}=1, 0.25$, centered at $x_{0}=-600,-500$ respectively. Tallest emerging wave. Time behaviour of: (a) Amplitude; (b) Speed.}
\label{fnls2_FIG10}
\end{figure}
\begin{figure}[htbp]
\centering
\subfigure[]
{\includegraphics[width=6.2cm]{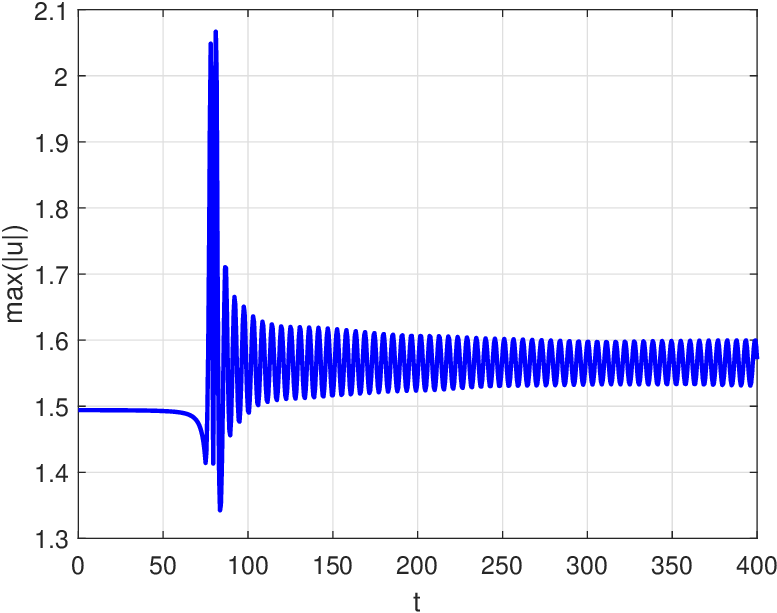}}
\subfigure[]
{\includegraphics[width=6.2cm]{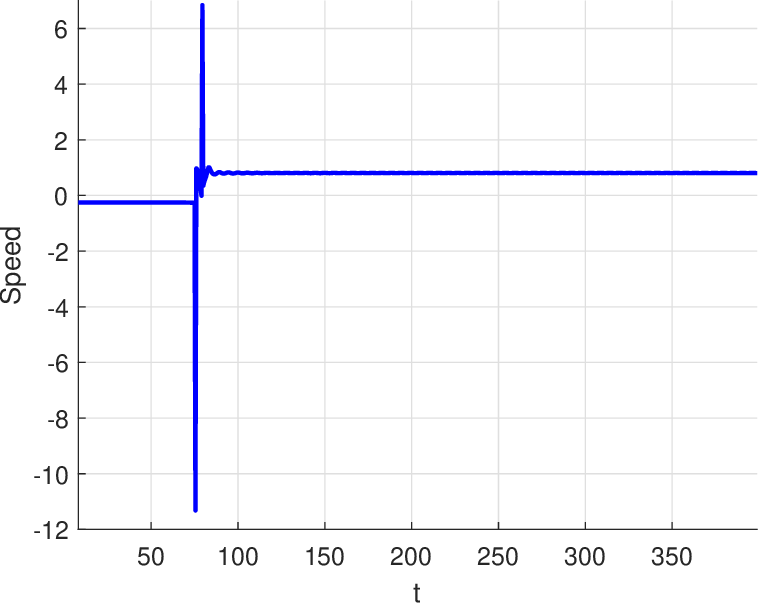}}
\caption{Head-on collisions. Numerical solution from a superposition of two computed solitary waves for $\theta(x)=A(x-x_{0})$, $a$ satisfying (\ref{aux1}), with $\lambda_{0}^{1}=1$ and $\lambda_{0}^{2}=1, -0.25$, centered at $x_{0}=-600,-500$ respectively. Tallest emerging wave. Time behaviour of: (a) Amplitude; (b) Speed.}
\label{fnls2_FIG12}
\end{figure}

The study of the interactions in the fractional case may be made here numerically. The first experiment is depicted in Figure \ref{fnls2_FIG9}. This represents the evolution of the modulus of the numerical approximation from an initial data given by the superposition of two computed solitary wave profiles, for $\theta(x)=A(x-x_{0})$, with $\lambda_{0}^{1}=1$ and $\lambda_{0}^{2}=1, 0.25$, centered at $x_{0}=-600,-500$, respectively, being the second one slower but larger in amplitude.  Thus, the waves interact at $t\approx 120$. After the collision, small-amplitude tails are generated behind and in front of each emerging wave, and affecting the corresponding amplitude and speed, cf, Figure \ref{fnls2_FIG10}, which displays the evolution of these parameters for the slower (and taller) wave.
Besides the inelastic character of the interaction (for example, the speed of the tallest wave decreases from $\lambda_{0}^{2}=0.25$ to $\lambda_{0}^{2}\sim 0.2421$), the behaviour of the amplitude in Figure \ref{fnls2_FIG10}(a) suggests that the slower wave may change of nature from the solitary wave to breather, due to the apparent time oscillatory behaviour.
\begin{figure}[htbp]
\centering
\subfigure
{\includegraphics[width=6.2cm]{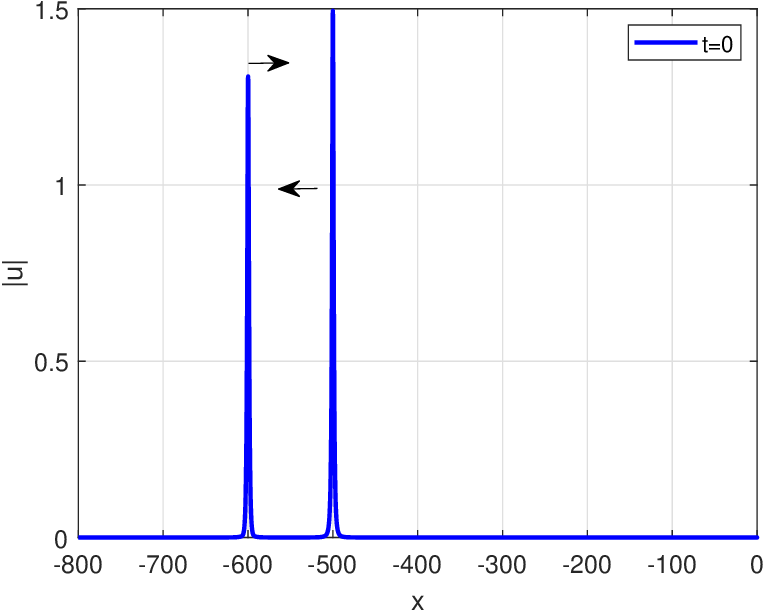}}
\subfigure
{\includegraphics[width=6.2cm]{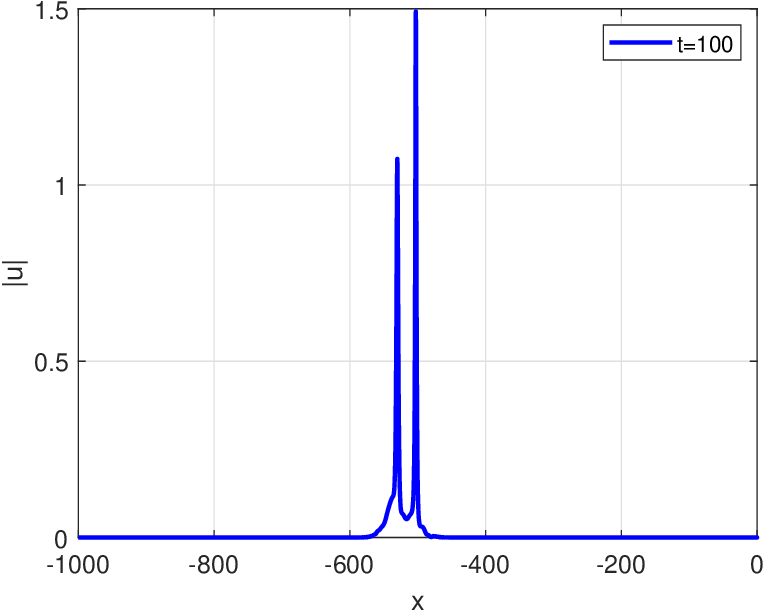}}
\subfigure
{\includegraphics[width=6.2cm]{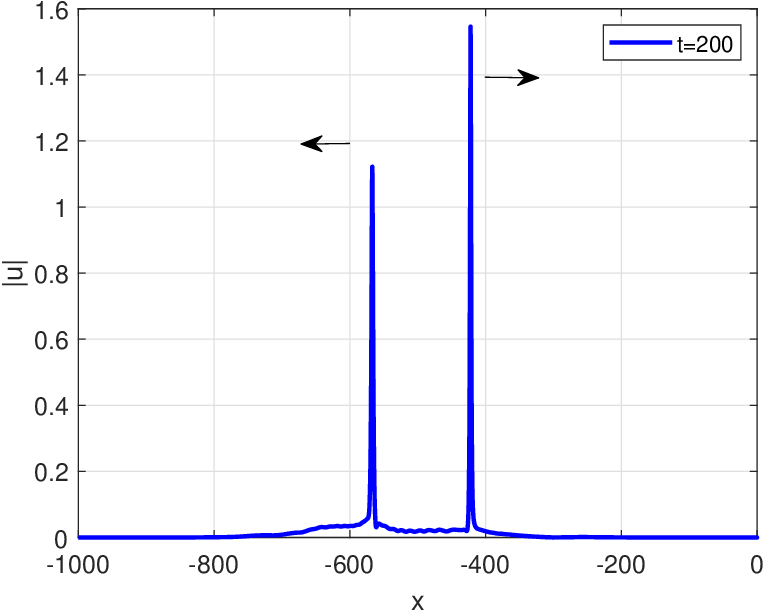}}
\subfigure
{\includegraphics[width=6.2cm]{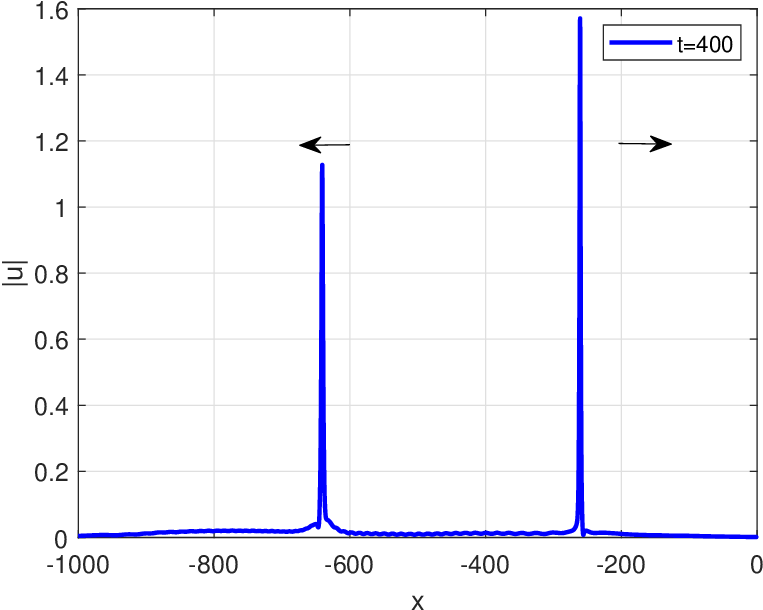}}
\caption{Head-on collisions. Evolution of the modulus of the numerical solution from a superposition of two computed solitary waves for $\theta(x)=A(x-x_{0})$, $a$ satisfying (\ref{aux1}), with $\lambda_{0}^{1}=1$ and $\lambda_{0}^{2}=1, -0.25$, centered at $x_{0}=-600,-500$ respectively.}
\label{fnls2_FIG11}
\end{figure}
Similar conclusions are suggested from experiments of head-on collisions. Note the existence of solitary wave solutions of (\ref{fnls1d}) holds for speeds $c_{s}=\lambda_{0}^{2}$ satisfying the condition (\ref{fnls_236b}). Then, we may consider as initial condition the superposition of two computed solitary wave profiles with speed $\lambda_{0}^{2}=-0.25$ (thus traveling to the left), centered at $x_{0}=-500$, and $\lambda_{0}^{2}=1$ (thus traveling to the right), centered at $x_{0}=-600$. The evolution of the numerical solution is monitored in Figure \ref{fnls2_FIG11}, with the nonlinear effects of the inelastic interaction on the emerging wave traveling to the left observed in the computation of the amplitude and speed as function of time in Figure \ref{fnls2_FIG12}. Note the change of the time behaviour in the amplitude and in the sign of the speed; after the collision, the waves seem to repel each other and travel in opposite directions.
\subsection{Complex interactions}
\label{sec43}
The robustness of the solitary waves can also be studied from their interactions with other types of waves. By way of illustration, a couple of examples from the experiments described in appendix \ref{appB} is mentioned here. They involve, respectively, collisions of solitary waves with large amplitude oscillatory waves (\ref{fnls_b2}) and interactions with waveforms by components (\ref{fnls_b4a}), (\ref{fnls_b4b}), that can be considered as nonsymmetric perturbations of the solitary waves when the parameters are small.

\begin{figure}[htbp]
\centering
\subfigure[]
{\includegraphics[width=6.2cm]{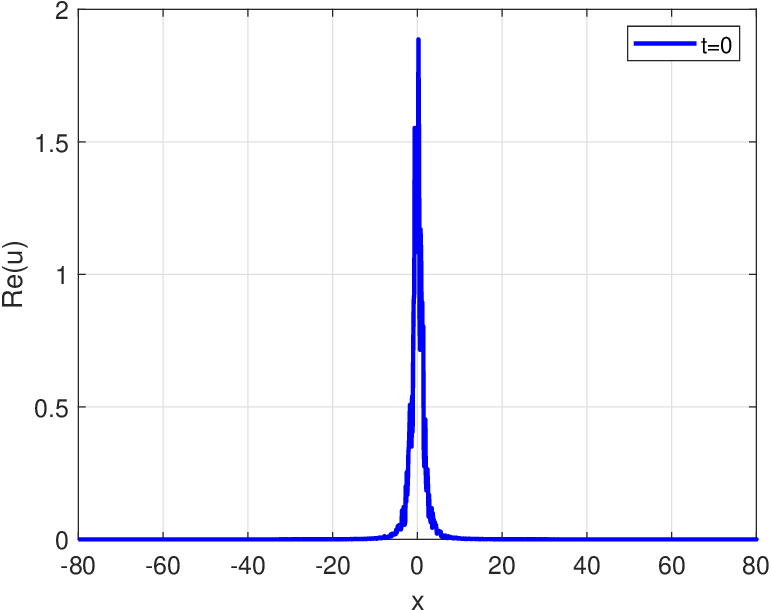}}
\subfigure[]
{\includegraphics[width=6.2cm]{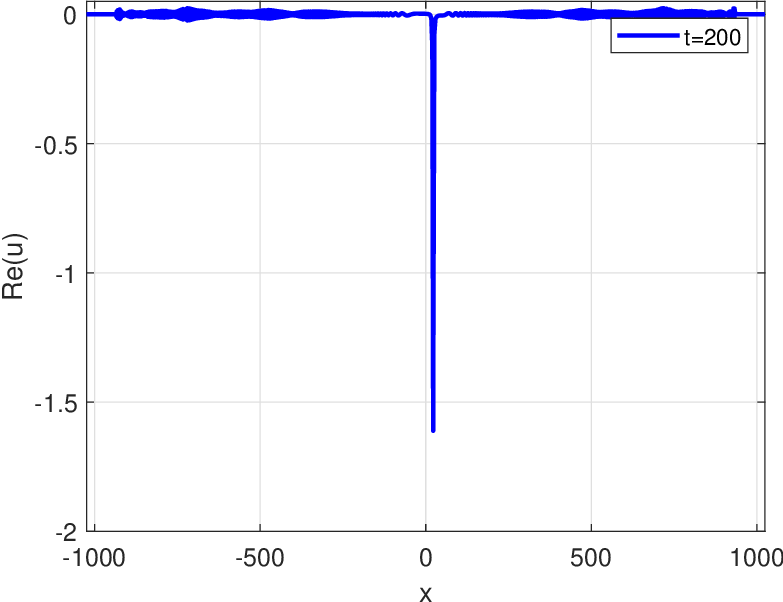}}
\subfigure[]
{\includegraphics[width=6.2cm]{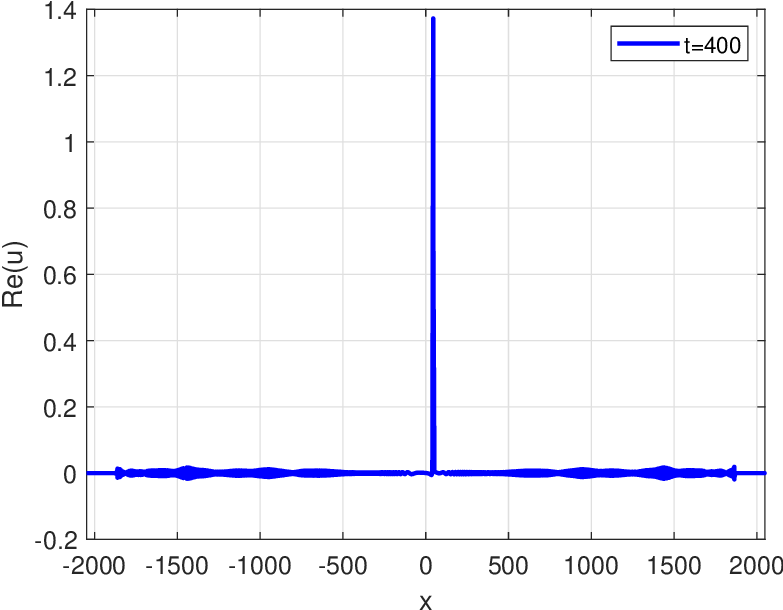}}
\subfigure[]
{\includegraphics[width=6.2cm]{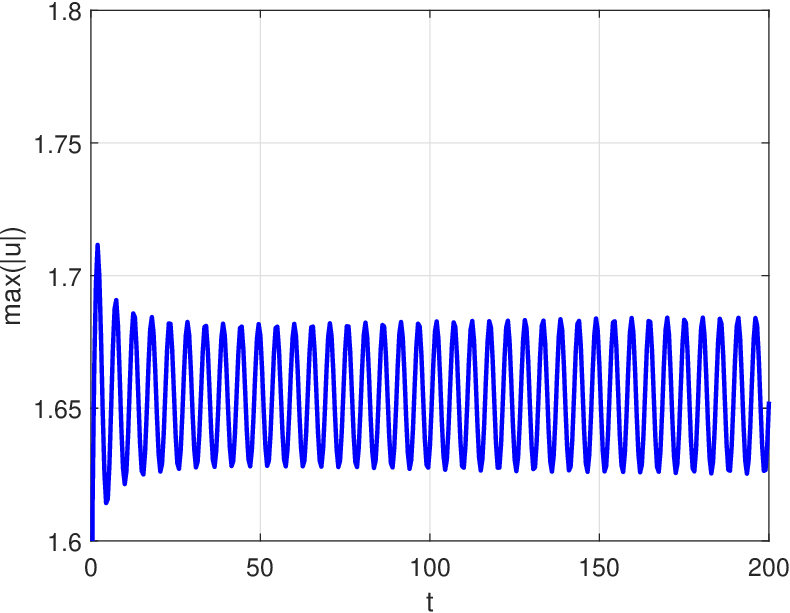}}
\caption{(a)-(c) Evolution of the $v$ component of the numerical solution from a perturbation of a solitary wave of the form  (\ref{fnls_b2}), (\ref{fnls_b3})  with $\beta=10^{7}$. (d) Time behaviour of the amplitude of the emerging wave.}
\label{fnls2_FIG13}
\end{figure}

\begin{figure}[htbp]
\centering
\subfigure[]
{\includegraphics[width=6.2cm]{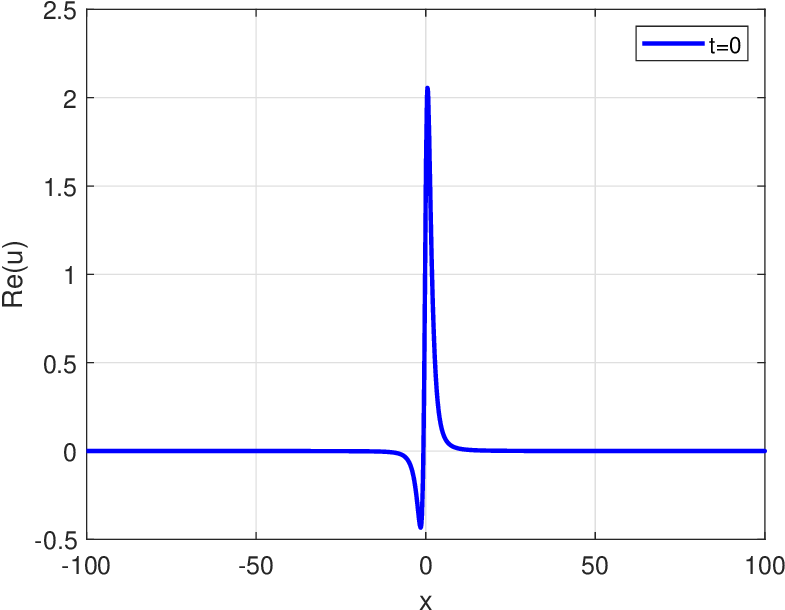}}
\subfigure[]
{\includegraphics[width=6.2cm]{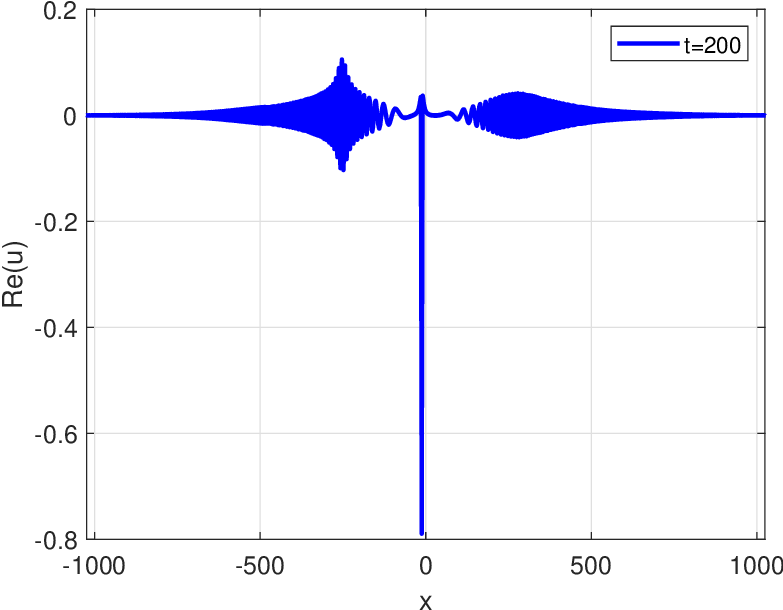}}
\subfigure[]
{\includegraphics[width=6.2cm]{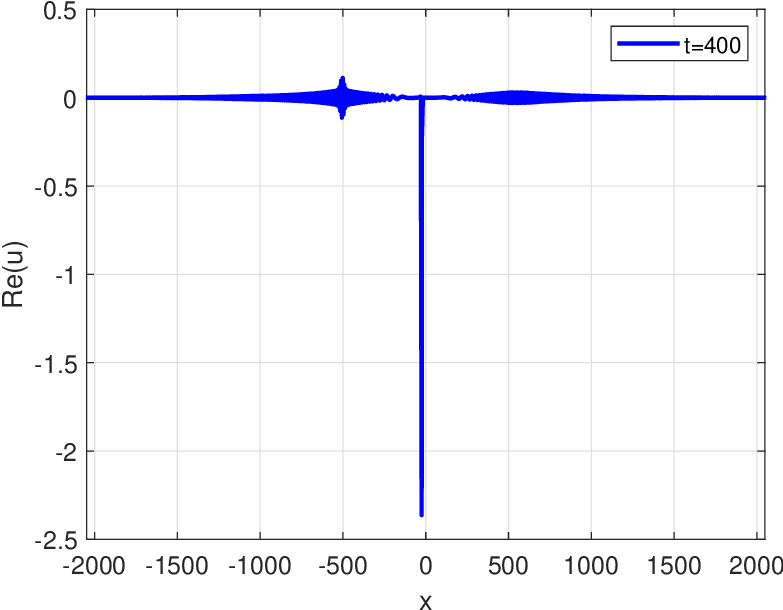}}
\subfigure[]
{\includegraphics[width=6.2cm]{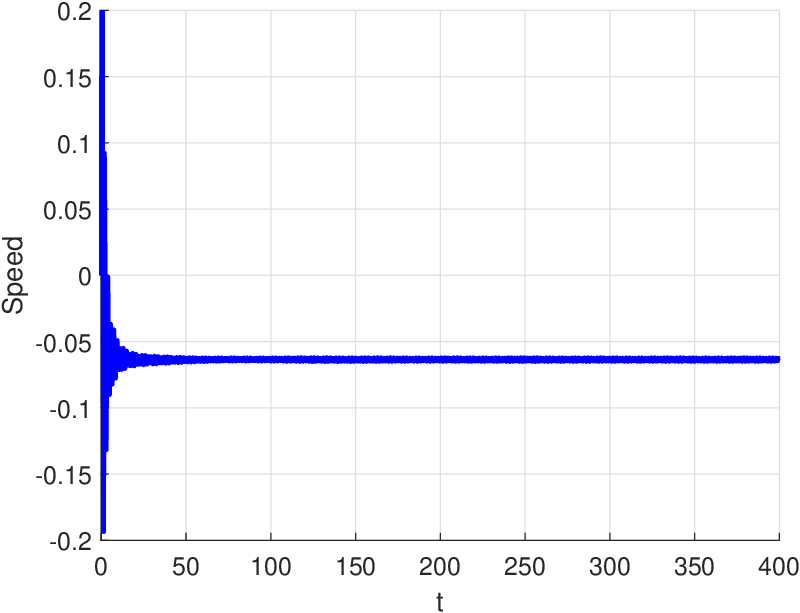}}
\caption{(a)-(c) Evolution of the $v$ component of the numerical solution from a perturbation of a solitary wave of the form  (\ref{fnls_b4a}), (\ref{fnls_b4b})  with $\alpha=3$. (d) Time behaviour of the speed of the emerging wave.}
\label{fnls2_FIG14}
\end{figure}

In the first experiment, the initial condition is a perturbed approximate solitary wave profile of the form  (\ref{fnls_b2}) with $\lambda_{0}^{1}=1, \lambda_{0}^{2}=0.25$, where a numerical noise function  (\ref{fnls_b3}) of relatively large amplitude (of order of $10^{-1}$, with $\beta=10^{7}$) is added to the $v$ component. The evolution of the corresponding component of the nuymerical solution is shown in Figure \ref{fnls2_FIG13}. The behaviour can be compared to that observed in section \ref{sec31} of perturbations of solitary wave profiles, with a larger amplitude of the emerging wave and some slower speed (from $\lambda_{0}^{2}=0.25$ of the initial profile to approximately $0.225$). 
The second experiment is illustrated in Figure \ref{fnls2_FIG14}, which shows the evolution of the $v$ component of the numerical solution from a perturbation of a solitary wave of the form  (\ref{fnls_b4a}), (\ref{fnls_b4b})  with $\alpha=3$. Note that the nonsymmetric perturbation generates a main wave which seems to have a breather-type behaviour and, as observed in Figure \ref{fnls2_FIG14}(d), travels to the left.
As a final comment, already outlined in the introduction, we made several experiments with large perturbations of different type, under the conditions $\sigma\geq 2s$ and initial data with negative energy (\ref{fnls3cc}), with the purpose of investigating numerically blow-up phenomena. In all the computations, the initial, perturbed wave evolved to the generation of some relatively large amplitude, moving breather, and no singularity ocurred.

\section*{Acknowledgments}
This research has been supported by Ministerio de Ciencia e Innovaci\'on project PID2023-147073NB-I00.

\bibliographystyle{plain}

\appendix
\section{The numerical method}
\label{appA}
The numerical study on the dynamics of solitary wave solutions of (\ref{fnls1d}) presented in sections \ref{sec2}-\ref{sec4} has been performed, via the method of lines, with a numerical scheme, for the periodic ivp (\ref{fnls1c}), based on a Fourier spectral method for the discretization in space and a fourth-order Runge-Kutta (RK) of composition type as time integrator. Stability and convergence properties of the full discretization have been established, either theoretically or numerically, when approximating other nonlinear dispersive models, cf. e.~g. \cite{DD2021} and references therein. In this appendix, the scheme is formulated and several experiments for checking the accuracy are shown.

The spatial discretization of  (\ref{fnls1c}),  (\ref{fnls1cc}) is a spectral Fourier Galerkin method on the space
\begin{eqnarray*}
S_{N}={\rm span}\{e^{i\widetilde{k}x}: \widetilde{k}=k\pi/L, k=-N,\ldots,N\},
\end{eqnarray*}
($N\geq 1$) of trigonometric polynomials on $[-L,L]$. If $T>0$, we approximate the solution $u=(v,w)$ of (\ref{fnls1c}),  (\ref{fnls1cc}) by real-valued functions $(v_{N},w_{N)}:[0,T]\rightarrow S_{N}\times S_{N}$ satisfying, for $0<t\leq T$ and $\varphi,\psi\in S_{N}$, the semidiscrete equations
\begin{eqnarray}
(v_{N_{t}},\varphi)-((-\partial_{xx})^{s}w_{N},\varphi)+((v_{N}^{2}+w_{N}^{2})^{\sigma}w_{N},\varphi)&=&0,\nonumber\\
-(w_{N_{t}},\psi)-((-\partial_{xx})^{s}v_{N},\psi)+((v_{N}^{2}+w_{N}^{2})^{\sigma}v_{N},\psi)&=&0,\label{fnlsA1}
\end{eqnarray}
and for $t=0$
\begin{eqnarray}
v_{N}(0)=P_{N}\widetilde{v}_{0},\quad w_{N}(0)=P_{N}\widetilde{w}_{0},\label{fnlsA2}
\end{eqnarray}
where $P_{N}$ denotes the $L^{2}$-projection operator on $S_{N}$ and
\begin{eqnarray*}
(\varphi,\psi)=\int_{-L}^{L}\varphi(x)\psi(x)dx,\quad \varphi,\psi\in L^{2}([-L,L]),
\end{eqnarray*}
is the $L^{2}$-inner product in $[-L,L]$.

The ivp (\ref{fnlsA1}), (\ref{fnlsA2}) can be formulated in terms of the Fourier coefficients 
\begin{eqnarray*}
\widehat{v}_{N}=\widehat{v}_{N}(k,t),\quad \widehat{w}_{N}=\widehat{w}_{N}(k,t),\quad -N\leq k\leq N, t\geq 0,
\end{eqnarray*}
of $v_{N}$ and $w_{N}$, as an ode ivp, for $ -N\leq k\leq N, 0<t\leq T$,
\begin{eqnarray}
\frac{d}{dt}\begin{pmatrix}\widehat{v}_{N}\\\widehat{w}_{N}\end{pmatrix}=F(\widehat{v}_{N},\widehat{w}_{N})&=&\begin{pmatrix}|\widetilde{k}|^{2s}\widehat{w}_{N}-\widehat{(v_{N}^{2}+w_{N}^{2})^{\sigma}w_{N}}\\
-|\widetilde{k}|^{2s}\widehat{v}_{N}+\widehat{(v_{N}^{2}+w_{N}^{2})^{\sigma}v_{N}}\end{pmatrix}, \label{fnlsA3}\\
\widehat{v}_{N}(k,0)=\widehat{\widetilde{v}_{0}}(k),&&
\widehat{w}_{N}(k,0)=\widehat{\widetilde{w}_{0}}(k),
\label{fnlsA4}
\end{eqnarray}
which is assumed to admit a unique solution, at least locally in time. The ode system (\ref{fnlsA3}), (\ref{fnlsA4}) is then numerically integrated in time with the $4$th-order, diagonally implicit Runge-Kutta composition method of the family of tableau
\begin{equation*}\label{44}
\begin{tabular}{c | c}
& $a_{ij}$\\ \hline
 & $b_{i}$
\end{tabular}=
\begin{tabular}{c | ccccc}
& $b_{1}/2$ & & &&\\
&  $b_{1}$ & $b_{2}/2$ &&&\\
&$b_{1}$&$b_{2}$&$\ddots$&&\\
&$\vdots$&$\vdots$ &&$\ddots$&\\
&$b_{1}$&$b_{2}$&$\cdots$&$\cdots$&$b_{s}/2$\\ \hline
 & $b_{1}$&$b_{2}$&$\cdots$&$\cdots$&$b_{s}$
\end{tabular},
\end{equation*}
in the particular case of $s=3$ stages, for which
\begin{eqnarray*}
&&b_{1}=(2+2^{1/3}+2^{-1/3})/3=\frac{1}{2-2^{1/3}}\sim 1.351,\nonumber\\
&& b_{2}=1-2b_{1}\sim -1.702,\quad b_{3}=b_{1}.\label{46}
\end{eqnarray*}
In addition, the scheme is of simple implementation using FFT techniques, and possesses several geometric and conservation properties, \cite{DR3}. The efficiency of the full discretization has been checked when it was used in the approximation of other nonlinear dispersive wave models, including stability issues. For the examples below, conditions of the form $N\Delta t\approx C$, $C$ constant (where $\Delta t$ denotes the time-step size), were checked to be enough for ensuring stability and convergence of the scheme.

\begin{figure}[htbp]
\centering
\subfigure
{\includegraphics[width=6.2cm]{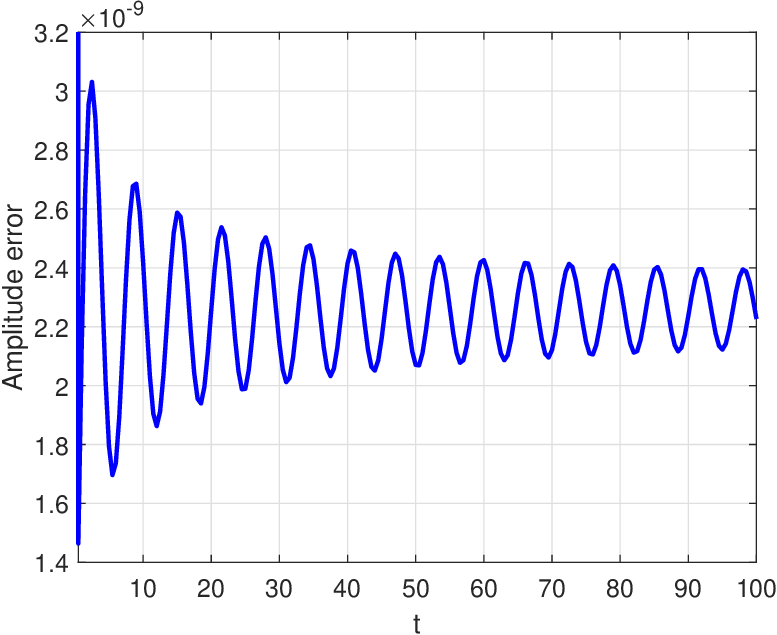}}
\subfigure
{\includegraphics[width=6.2cm]{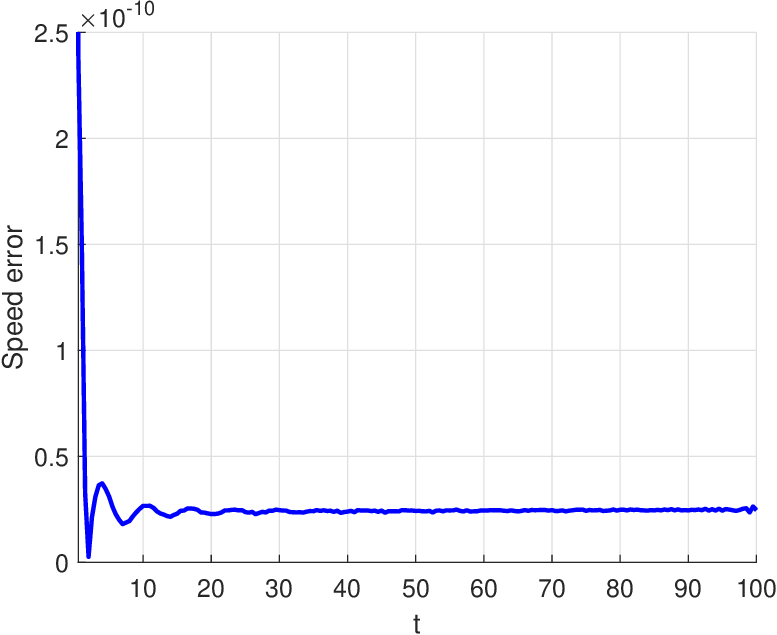}}
\caption{Time behaviour of the errrors in (a) amplitude; (b) speed w.r.t. the solitary-wave solution (\ref{fnls_A5}) with $\sigma=1, \lambda_{0}^{1}=1, \lambda_{0}^{2}=0.25$, with $\Delta t=6.25\times 10^{-3}$.}
\label{fnls2_FIG00}
\end{figure}

We now present some numerical experiments to validate the performance of the fully discrete method introduced above  and to give some confidence on the accuracy of the computations made in sections \ref{sec2}-\ref{sec4}. They were made to approximate the classical NLS, $s=1$ in (\ref{fnls1d}), and its solitary wave solutions of the form (\ref{fnls22_7}) where, \cite{DuranS2000}
\begin{eqnarray}
\rho(x)&=&(a(\sigma +1))^{1/2\sigma}({\rm sech}\sigma\sqrt{a}x)^{1/\sigma},\quad a=\lambda_{0}^{1}-\frac{(\lambda_{0}^{2})^{2}}{4},\nonumber\\
\theta(x)&=&\frac{\lambda_{0}^{2}}{2}x.\label{fnls_A5}
\end{eqnarray}
The numerical solution at $T=100$ is first compared with the exact solution (\ref{fnls22_7}), (\ref{fnls_A5}) with $\sigma=1, \lambda_{0}^{1}=1, \lambda_{0}^{2}=0.25, x_{0}=\theta_{0}=0$. The errors in the $L^{2}$ norm of the $v$ and $w$ components, for several time-step sizes are displayed in Table \ref{texpe1}, showing the fourth order of convergence of the time discretization. (Note that, since the soliton solutions are smooth, a spectral order of convergence of the semidiscrete approximation is expected, \cite{DD2021}.)
\begin{table}
\begin{tabular}{c|c|c|c|c|}
$\Delta t$&$v$ Error&Rate&$w$ Error&Rate\\
\hline
$2.5\times 10^{-2}$&$1.1621\times 10^{-4}$&&$1.9446\times 10^{-4}$&\\
$1.25\times 10^{-2}$&$7.2820\times 10^{-6}$&$3.9963$&$1.2187\times 10^{-5}$&$3.9961$\\
$6.25\times 10^{-3}$&$4.5630\times 10^{-7}$&$3.9963$&$7.6359\times 10^{-7}$&$3.9963$\\
$3.125\times 10^{-3}$&$2.7478\times 10^{-8}$&$4.0537$&$4.6003\times 10^{-8}$&$4.0530$\\
\hline
\end{tabular}
\caption{$L^{2}$ errors and temporal convergence rates. Solitary-wave solution (\ref{fnls_A5}) with $\sigma=1, \lambda_{0}^{1}=1, \lambda_{0}^{2}=0.25$,  $T=100$, $N=4096$.\label{texpe1}}
\end{table}
The accuracy is also observed when simulating relevant parameters of the soliton, such as the amplitude and speed. The time behaviour of the errors with $\Delta t=6.25\times 10^{-3}$ is depicted in Figure \ref{fnls2_FIG00}. Other experiments, concerning the preservation of discrete versions of the invariants (\ref{fnls3aa})-(\ref{fnls3cc}) (not shown here), may illustrate the benefits of the geometric properties of the scheme, cf. \cite{DR1}.

\section{Types of experiments}
\label{appB}
In this appendix we make a description of the experiments performed in the present paper. All the experiments share the following values for the parameters
$$N=16384, L=1024, T=400, s=0.8, \sigma=1, \lambda_{0}^{1}=1,$$ although other similar experiments, not shown here, were performed, using different values of $s$ and $\sigma$, to confirm a similar behaviour of the output. We also consider different time-step sizes $\Delta t$ and make use of two different types of phases
\begin{eqnarray*}
\theta(x)=A(x-x_{0}),\quad \theta(x)=(x-x_{0})^{2},
\end{eqnarray*}
with $A$ satisfying (\ref{aux1}), in order to check a potential influence of the properties of the solitary wave profiles in their dynamics, \cite{DR1}. The full set of experiments is included in a technical report available from the authors upon request.

The choice of the perturbations is inspired in previous studies, \cite{DDLM}. They are the following.
\subsection{Perturbations in amplitude}
We take $\lambda_{0}^{2}=0.25$ or $0.5$, $x_{0}=\theta_{0}=0$, and generate numerically an approximate solitary wave profile $(\widetilde{v},\widetilde{w})$ by using the procedure described in \cite{DR1}. Then the system (\ref{fnls1b}) is integrated with initial conditions of the form
\begin{eqnarray}
\widetilde{v}_{0}(x)=A_{1}v_{0}(x),\quad \widetilde{w}_{0}(x)=A_{2}w_{0}(x),\label{fnls_b1}
\end{eqnarray}
for several values of the perturbation factors $A_{j}, j=1,2$.
\subsection{Perturbations with a numerical noise}
\begin{figure}[htbp]
\centering
\subfigure[]
{\includegraphics[width=6.2cm]{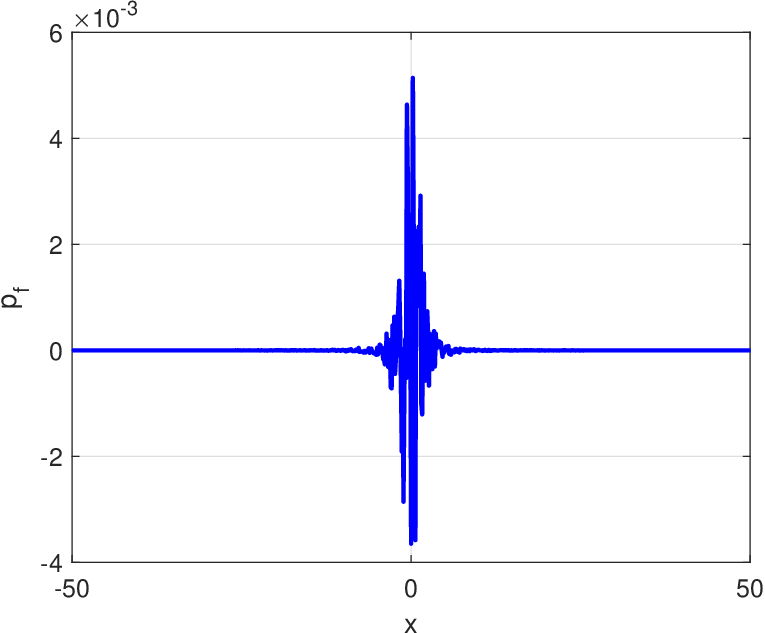}}
\subfigure[]
{\includegraphics[width=6.2cm]{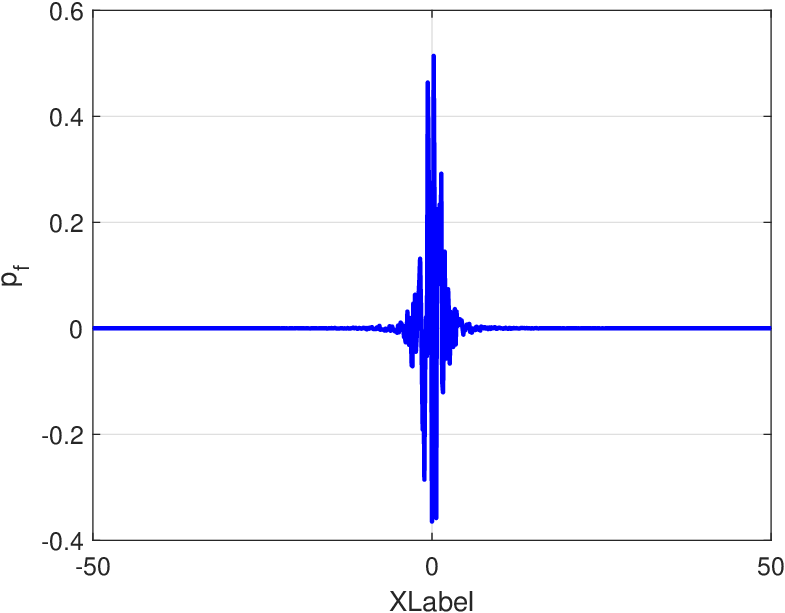}}
\caption{Numerical noise  (\ref{fnls_b3}). (a) $m=10^{5}$; (b) $m=10^{7}$.}
\label{fnlsB1}
\end{figure}

The initial conditions in this case are of the form
\begin{eqnarray}
&&v_{s}=\widetilde{v}+p_{v},\quad w_{s}=\widetilde{w},\; {\rm or}\nonumber\\
&&v_{s}=\widetilde{v},\quad w_{s}=\widetilde{w}+p_{w},\label{fnls_b2}
\end{eqnarray}
where $p_{f}, f=v,w$, in (\ref{fnls_b2}) represents a numerical noise function which depends on $u$ as
\begin{eqnarray}
p_{f}=p_{f}(x,m)=m({\rm double}(f(x,0))-{\rm single}(f(x,0))),\label{fnls_b3}
\end{eqnarray}
where ${\rm double}(f(x,0))$ is the double precision function $f(x,0)$, ${\rm single}(f(x,0))$ is the single precision function $f(x,0)$, and $m$ is a parameter determining the size of the noise. Figure \ref{fnlsB1} represents the function (\ref{fnls_b3}) with $m=10^{5}, 10^{7}$.
\subsection{Nonsymmetric perturbations}
The initial conditions in this case are of the form
\begin{eqnarray}
&&v_{s}=\widetilde{v}p,\quad w_{s}=\widetilde{w},\; {\rm or}\nonumber\\
&&v_{s}=\widetilde{v},\quad w_{s}=\widetilde{w}p,\label{fnls_b4a}
\end{eqnarray}
where
\begin{eqnarray}
p(x)=1+\alpha{\rm tanh}\left(\frac{1}{2}(x-x_{0})\right),\quad \alpha\in\mathbb{R}.\label{fnls_b4b}
\end{eqnarray}
\subsection{Experiments for resolution property}
The resolution property is illustrated with two types of initial conditions:
\begin{itemize}
\item[(1)] Of the form (\ref{fnls_b1}) for large values of $(A_{1},A_{2})$.
\item[(2)] Initial data of modulated Gaussian type
\begin{eqnarray}
&&\rho(x)=A_{1}e^{-A_{2}(x-x_{0})^{2}},\nonumber\\
&&v_{s}(x)=\rho(x)\cos(\theta(x)+\theta_{0}),\quad w_{s}(x)=\rho(x)\sin(\theta(x)+\theta_{0}),\label{fnls_b5}
\end{eqnarray}
and several values of $A_{j}, j=1,2$.
\end{itemize}
\subsection{Perturbations of the parameter $s$}
Several experiments were also performed from an initial approximate solitary wave profile $(\widetilde{v}_{\overline{s}},\widetilde{w}_{\overline{s}})$, generated with a value $\overline{s}\in (1/2,1)$ of the parameter of the fractional Laplacian. Then the parameter $s$ of (\ref{fnls1b}) is perturbed as $s=\overline{s}+\epsilon, \epsilon>0$. The corresponding numerical approximation is evaluating the dynamics of solitary waves of equations (\ref{fnls1b}) which are close. The experiments, not shown in the present paper, are included in the technical report available from the authors upon request.
\subsection{Experiments of interactions}
The experiments developed in sections \ref{sec42} and \ref{sec43} involve interactions between approximate solitary waves (overtaking and head-on) and interactions between one approximate solitary wave and other types of waveforms. In the first case, the initial condition is a superposition of computed solitary waves of different speed; in the second one, the initial data consists of a superposition of an approximate solitary wave plus a wave profile of the form of a numerical noise function defined in (\ref{fnls_b3}) or the nonsymmetric form given by (\ref{fnls_b4a}), (\ref{fnls_b4b}).

\end{document}